\documentclass[11pt]{article}
\usepackage[authoryear,round]{natbib}
\usepackage{graphicx}
\usepackage{latexsym}
\usepackage[left=1.5in,top=1.5in,right=1.5in,bottom=1.5in]{geometry}
\usepackage{amsmath}
\usepackage{amssymb}
\usepackage{booktabs}
\usepackage{color}
\usepackage{tikz}
\usepackage{import}
\usepackage[ruled,linesnumbered]{algorithm2e}
\usepackage{savesym}
\usepackage{mdframed}
\mdtheorem{egbox}{Example}
\newtheorem{remark}{Remark}[section]

\usetikzlibrary{shapes}
\usetikzlibrary{arrows}
\definecolor{darkblue}{rgb}{0,0.4,0.9}
\definecolor{gray10}{rgb}{0.1,0.1,0.1}
\definecolor{gray20}{rgb}{0.2,0.2,0.2}
\definecolor{gray30}{rgb}{0.3,0.3,0.3}
\definecolor{gray40}{rgb}{0.4,0.4,0.4}
\definecolor{gray60}{rgb}{0.6,0.6,0.6}
\definecolor{gray80}{rgb}{0.8,0.8,0.8}
\definecolor{gray90}{rgb}{0.9,0.9,.9}
\definecolor{gray95}{rgb}{0.95,0.95,.95}
\definecolor{gray96}{rgb}{0.96,0.96,.96}
\definecolor{lgreen} {RGB}{180,210,100}
\definecolor{dblue}  {RGB}{20,66,129}
\definecolor{ddblue} {RGB}{11,36,69}
\definecolor{lred}   {RGB}{220,0,0}
\definecolor{nred}   {RGB}{224,0,0}
\definecolor{norange}{RGB}{230,120,20}
\definecolor{nyellow}{RGB}{255,221,0}
\definecolor{ngreen} {RGB}{98,158,31}
\definecolor{dgreen} {RGB}{78,138,21}
\definecolor{nblue}  {RGB}{28,130,185}
\definecolor{jblue}  {RGB}{20,50,100}
\definecolor{nnyellow}{RGB}{235,200,0}
\definecolor{purple}{RGB}{150, 0, 120}
\definecolor{sgGreen} {RGB}{20, 180, 50}
\definecolor{revised}{rgb}{0,0,0.9}
\newtheorem{definition}{Definition}

\newtheorem{theorem}{Theorem}
\newtheorem{corollary}{Corollary}
\newtheorem{lemma}{Lemma}

\savesymbol{nl}

\restoresymbol{algo}{nl}
\newcommand{\pl}{\|}
\newcommand{\openr}{\hbox{${\rm I\kern-.2em R}$}}
\newcommand{\openn}{\hbox{${\rm I\kern-.2em N}$}}
\DeclareMathOperator*{\argmax}{arg\,max}
\DeclareMathOperator*{\argmin}{arg\,min}

\title{Nonparametric Bootstrap Inference for the Targeted Highly Adaptive LASSO Estimator}

\author{Weixin Cai and Mark van der Laan
 \\ Division of Biostatistics, University of California, Berkeley\\ {\tt laan@berkeley.edu}
\\
} \date{\today}

\begin{document}
\maketitle
 
\begin{abstract}
The Highly-Adaptive-LASSO Targeted Minimum Loss Estimator (HAL-TMLE) is an efficient plug-in estimator of a pathwise differentiable parameter in a statistical model that at minimal (and possibly only) assumes that the sectional variation norm of the true nuisance functions (i.e., relevant part of data distribution) are finite. It relies on an initial estimator (HAL-MLE) of the nuisance functions by minimizing the empirical risk over the parameter space under the constraint that the sectional variation norm of the candidate functions are bounded by a constant, where this constant can be selected with cross-validation. In this article we establish that the nonparametric bootstrap for the HAL-TMLE, fixing the value of the sectional variation norm at a value larger or equal than the cross-validation selector, provides a consistent method for estimating the normal limit distribution of the HAL-TMLE. 

In order to optimize the finite sample coverage of the nonparametric bootstrap confidence intervals, we propose a selection method for this sectional variation norm that is based on running the nonparametric bootstrap for all values of the sectional variation norm larger than the one selected by cross-validation, and subsequently determining a value at which the width of the resulting confidence intervals reaches a plateau.  

We demonstrate our method for 1) nonparametric estimation of the average treatment
effect when observing a covariate vector, binary treatment,
and outcome, and for 2) nonparametric estimation of the integral of the square of the multivariate density of the data distribution. In addition, we also present simulation results for these two examples demonstrating the excellent finite sample coverage of bootstrap-based confidence intervals.

\end{abstract}


\section{Introduction}

We consider estimation of a pathwise differentiable real valued target estimand based on observing $n$ independent and identically distributed observations $O_1,\ldots,O_n$ from a data distribution $P_0$ known to belong to a statistical model ${\cal M}$. 
A target parameter $\Psi:{\cal M}\rightarrow \openr$ is a  mapping that maps a possible data distribution $P\in {\cal M}$ into real number, while $\psi_0=\Psi(P_0)$ represents the statistical estimand. 
The canonical gradient $D^*(P)$ of the pathwise derivative of the target parameter at a distribution $P$ defines an asymptotically efficient estimator among the class of regular estimators \citep{Bickeletal97}: an estimator $\psi_n$ is asymptotically efficient at $P_0$ if and only if it is asymptotically linear at $P_0$ with influence curve $D^*(P_0)$:
\[
\psi_n-\psi_0=\frac{1}{n}\sum_{i=1}^n D^*(P_0)(O_i)+o_P(n^{-1/2}).\]
The target parameter $\Psi(P)$ depends on the data distribution $P$ through a parameter $Q=Q(P)$, while the canonical gradient $D^*(P)$ possibly also depends on another nuisance parameter $G(P)$: $D^*(P)=D^*(Q(P),G(P))$. Both of these nuisance parameters are chosen so that they can be defined as a minimizer of the expectation of a specific loss function: $Q(P))=\argmin_{Q\in Q({\cal M})}PL_1(Q)$ and $G(P)=\argmin_{G\in G({\cal M})}P L_2(G)$, where we used the notation $Pf\equiv \int f(o)dP(o)$. We consider the case that the parameter spaces $Q({\cal M})=\{Q(P):P\in {\cal M}\}$ and $G({\cal M})=\{G(P):P\in {\cal M}\}$  for these nuisance parameters $Q$ and $G$ are contained in the set of multivariate cadlag functions with sectional variation norm $\pl \cdot\pl_v^*$ \citep{Gill&vanderLaan&Wellner95} bounded by a constant (this norm will be defined in the next section). 

We consider a targeted minimum loss-based (substitution) estimator $\Psi(Q_n^*)$  \citep{vanderLaan&Rubin06,vanderLaan08, vanderLaan&Rose11,vanderLaan&Rose17}
of the target parameter that uses as initial estimator   of these nuisance parameters $(Q_0,G_0)$ the highly adaptive lasso minimum loss-based estimators (HAL-MLE)  $(Q_n,G_n)$ defined by minimizing the empirical mean of the loss over the parameter space \citep{vanderLaan15,Benkeser&vanderLaan16}.  Since the HAL-MLEs converge at a rate faster than $n^{-1/2}$ with respect to (w.r.t.) the loss-based quadratic dissimilarities (to be defined later, which corresponds with a rate faster than $n^{-1/4}$ for estimation of $Q_0$ and $G_0$), this HAL-TMLE has been shown to be asymptotically efficient under weak regularity conditions \citep{vanderLaan15}. Statistical inference could therefore be based on the normal limit distribution in which the asymptotic variance is estimated with an estimator of the variance of the canonical gradient. 
In that case, inference is ignoring the potentially very large contributions of the higher order remainder which could, in finite samples, easily dominate the first order empirical mean of the efficient influence curve term when the size of the nuisance parameter spaces is large (e.g., dimension of data is large and model is nonparametric). 

In this article we propose the nonparametric bootstrap to  obtain a better estimate of the  finite sample distribution of the HAL-TMLE than the normal limit distribution. The bootstrap  fixes the sectional variation norm at the values used for the HAL-MLEs $(Q_n,G_n)$ on a bootstrap sample. We propose a data adaptive selector of this tuning parameter tailored to obtain improved finite sample coverage for the resulting confidence intervals. 
\subsection{Organization}

In Section \ref{secttwo} we formulate the estimation problem and motivate the challenge for statistical inference. 
In Section \ref{sectthree} we  present the nonparametric bootstrap estimator of the actual sampling distribution of the HAL-TMLE which  thus incorporates estimation of its higher order stochastic behavior, and can thereby be expected to outperform the Wald-type confidence intervals.   We prove that this nonparametric bootstrap is asymptotically consistent for the optimal normal limit distribution. Our results also prove that the nonparametric bootstrap preserves the asymptotic behavior of the HAL-MLEs of our nuisance parameters $Q$ and $G$, providing further evidence for good performance of the  nonparametric bootstrap.
 Importantly, our results demonstrate that the approximation error of the nonparametric bootstrap estimate of the true finite sample distribution of the HAL-TMLE is mainly driven by the approximation  error of the nonparametric bootstrap for estimating the finite sample distribution of a well behaved empirical process. 
 In Section \ref{sectfour} we present a plateau selection method for selecting the fixed sectional variation norm in the nonparametric bootstrap and a bias-correction in order to obtain improved finite sample coverage for the resulting confidence intervals. 

In Section \ref{sectfive} we demonstrate our methods for two examples involving a nonparametric model and a specified target parameter: average treatment effect and integral of the square of the data density. 
In Section \ref{sectsix} we carry out a simulation study to demonstrate the practical performance of our proposed nonparametric bootstrap based confidence intervals w.r.t. their finite sample coverage.
We conclude with a discussion in Section \ref{sectdisc}. 
Proofs of our Lemma and Theorems  have been deferred to the Appendix.
We refer to our accompanying technical report for additional bootstrap methods and results based on applying the nonparametric bootstrap to an exact second order expansion of the HAL-TMLE, and to various upper bounds of this exact second order expansion.

\section{General formulation of statistical estimation problem and motivation for finite sample inference}\label{secttwo}

\subsection{Statistical model and target parameter}
Let $O_1,\ldots,O_n$ be $n$ i.i.d. copies of a random variable $O\sim P_0\in {\cal M}$. 
Let $P_n$ be the empirical probability measure of $O_1,\ldots,O_n$.
Let $\Psi:{\cal M}\rightarrow\openr$ be a real valued parameter that is pathwise differentiable at each $P\in {\cal M}$ with 
canonical gradient $D^*(P)$.  That is, given a collection of one dimensional  submodels $\{P_{\epsilon}^S:\epsilon\}\subset {\cal M}$ through $P$ at $\epsilon =0$ with score $S$, for each of these submodels the derivative $\left . \frac{d}{d\epsilon}\Psi(P_{\epsilon}^S)\right |_{\epsilon =0}$ can be represented as a covariance  $ E_P D(P)(O)S(O)$ of a gradient $D(P)$ with the score $S$.
The latter is  an inner product of a gradient $D(P)\in L^2_0(P)$ with the score $S$ in the Hilbert space $L^2_0(P)$ of functions of $O$ with mean zero (under $P$) endowed with inner product $\langle S_1,S_2\rangle_P=P S_2S_2$. Let $\pl f\pl_P\equiv \sqrt{\int f(o)^2dP(o)}$ be the Hilbert space norm.
Such an element $D(P)\in L^2_0(P)$ is called a gradient of the pathwise derivative of $\Psi$ at $P$.
The canonical gradient  $D^*(P)$ is the unique gradient that is an element of the tangent space defined as the closure of the linear span of the collection of scores generated by this family of submodels. 
Define the exact second-order remainder
\begin{equation}\label{exactremainder}
R_2(P,P_0)\equiv\Psi(P)-\Psi(P_0)+(P-P_0)D^*(P),\end{equation}
where $(P-P_0)D^*(P)=-P_0D^*(P)$ since $D^*(P)$ has mean zero under $P$.

\begin{egbox}[(Treatment-specific mean)]
    Let $O = (W, A, Y) \sim P_0 \in \mathcal{M}$, where $A \in \{0,1\}$ is a binary
    treatment, $Y \in\{0,1\}$ is a binary outcome, and ${\cal M}$ is a nonparametric model. For a possible data
    distribution $P$, let $\bar{Q}(P) = \mathbb{E}_P(Y|A,W)$ be the outcome regression, $G(P) = P(A = 1|W)$ be the propensity score, and let $Q_W(P)$ be the probability distribution of $W$. The treatment-specific mean parameter is defined by $\Psi(P) = \mathbb{E}_P\mathbb{E}_P(Y|A = 1, W)$. Let $Q=(\bar{Q},Q_W)$ and note that the data distribution $P$ is determined by $(Q,G)$. The canonical gradient of $\Psi$ at $P$ is $$D^*(P)=D^*(Q,G) =
    \frac{I(A=1)}{G(A|W)}(Y-\bar{Q}(A,W)) + \bar{Q}(1,W) - \Psi(Q).$$ The second-order remainder
    $R_2(P,P_0) \equiv \Psi(P) - \Psi(P_0) + P_0D^*(P)$ is given by:
    {
    \begin{align*}
    R_2(Q,G,Q_0,G_0) = \int {\frac{(G - G_0)(w)}{G(w)} (\bar{Q} -\bar{ Q}_0)(1,w) dP_0(w)}
    \end{align*}
    }
\end{egbox}

Let $Q:{\cal M}\rightarrow Q({\cal M})$ be a function valued parameter so that $\Psi(P)=\Psi_1(Q(P))$ for some $\Psi_1$.
For notational convenience, we will abuse  notation  by referring to the target parameter with $\Psi(Q)$ and $\Psi(P)$ interchangeably. 
Let $G:{\cal M}\rightarrow G({\cal M})$ be a function valued parameter so that $D^*(P)=D^*_1(Q(P),G(P))$ for some $D^*_1$. Again, we will use the notation $D^*(P)$ and $D^*(Q,G)$ interchangeably. 

For each $Q\in Q({\cal M})$, let $L_1(Q)$ be a function of $O$ so that
\begin{align*}
Q_0 = \argmin_{Q\in Q(\cal M)} P_0 L_1(Q).
\end{align*}
Similarly, for each $G\in G({\cal M})$, let $L_2(G)$ be a function of $O$ so that
\begin{align*}
G_0 = \argmin_{G\in G(\cal M)} P_0L_2(G).    
\end{align*}
We refer to $L_1(Q)$ and $L_2(G)$ as loss functions for $Q_0$ and $G_0$.  Let $d_{01}(Q,Q_0)=P_0L_1(Q)-P_0L_1(Q_0)\geq 0$ and $d_{02}(G,G_0)=P_0L_2(G)-P_0L_2(G_0)\geq 0$ be the loss-based dissimilarities for these two nuisance functions. The  loss based dissimilarity is often called the regret. 
 Assume that the loss functions are uniformly bounded in the sense that $\sup_{Q\in Q({\cal M}),O}| L_1(Q)(O)| <\infty$ and $\sup_{G \in G({\cal M}), O}| L_2(G)(O)| <\infty$.
In addition, assume
\begin{eqnarray}
\sup_{Q\in Q({\cal M})}\frac{P_0\{L_1(Q)-L_1(Q_0)\}^2}{d_{01}(Q,Q_0)}&<&\infty\nonumber \\
\sup_{G\in G({\cal M})}\frac{P_0\{L_2(G)-L_2(G_0)\}^2}{d_{02}(G,G_0)}&<&\infty
\label{M20bounds}.
\end{eqnarray}
This condition holds for most common bounded loss functions (such as mean-squared error loss and cross entropy loss), and it guarantees that the loss-based dissimilarities $d_{01}(Q,Q_0)$ and $d_{02}(G,G_0)$ behave as a square of an $L^2(P_0)$-norm. 
These two universal bounds on the loss function yield the oracle inequality for the cross-validation selector among a set of candidate estimators \citep{vanderLaan&Dudoit03,vanderVaart&Dudoit&vanderLaan06,vanderLaan&Dudoit&vanderVaart06,vanderLaan&Polley&Hubbard07,Chpt3}. In particular, it establishes that the cross-validation selector is asymptotically equivalent to the oracle selector. 

\begin{egbox}[(Treatment-specific mean)]
    For the treatment-specific mean parameter, the $\bar{Q}$ function is the outcome regression $E(Y|A,W)$, and $G = P(A=1|W)$ is the propensity score. The other component $Q_W$ of $Q$ will be estimated with the empirical probability measure, which is an NPMLE, so that a TMLE will not update this estimator. Let $L_1(\bar{Q})(O)=-\{Y\log \bar{Q}(A,W)+(1-Y)\log(1-\bar{Q}(A,W))\}$ be the negative log-likelihood loss for the outcome regression. Similarly, $L_2(G)$ is the negative-log-likelihood loss for propensity score. When, for some $\delta>0$, $G>\delta>0$ and $\delta<\bar{Q}<1-\delta$,  then the loss functions are uniformly bounded with finite universal bounds (\ref{M20bounds}).
\end{egbox}

\subsubsection{Donsker class condition}
Our formal theorems need to assume that $\{L_1(Q):Q\in Q({\cal M})\}$, $\{L_2(G):G\in G({\cal M})\}$, and $\{D^*(Q,G):Q\in Q({\cal M}), G\in G({\cal M})\}$ are uniform (in $P\in {\cal M}$) Donsker classes, or, equivalently, that the union ${\cal F}$ of these classes is a uniform Donsker class. We remind the reader that a covering number $N(\epsilon,{\cal F},L^2(\Lambda))$  is defined as the minimal number of balls of size $\epsilon$ w.r.t. $L^2(\Lambda)$-norm that are needed to cover the set ${\cal F}$ of functions embedded in $L^2(\Lambda)$.
Let $\alpha\in (0,1)$ be defined such that
\begin{equation}\label{alphad}
\sup_{\Lambda}\log^{1/2}(N(\epsilon,{\cal F},L^2(\Lambda)) = O( \epsilon^{-(1-\alpha)}).\end{equation}
Our formal results will refer to a rate of convergence of the HAL-MLEs w.r.t. loss based dissimilarity given by  $n^{-1/2-\alpha/4}$ implied by  this index $\alpha$ \citep{vanderLaan15}. 
In this article we will focus here on the following special Donsker class, in which case $\alpha$ can be chosen as $2/(d+2)$.

\subsubsection{Loss functions and canonical gradient have a uniformly bounded sectional variation norm}
We assume that the loss functions and canonical gradient are cadlag functions with a universal bound on the sectional variation norm. The latter class of functions is indeed a uniform Donsker class. In the sequel we will assume this, but we remark here that throughout we could have replaced this class of cadlag functions with a universal bound on the sectional variation norm by any other uniform Donsker class.  Below we will present a particular class of models ${\cal M}$ in which we assume that the nuisance parameters $Q$ and $G$  themselves fall in such classes of functions, so that generally also $L_1(Q), L_2(G)$ and $D^*(Q,G)$ will fall in this class. All our applications have been covered by the latter type of models.

We will formalize this condition now.
Suppose that $O\in [0,\tau]\subset\openr^d_{\geq 0}$ is a $d$-variate random variable with support contained in a $d$-dimensional cube $[0,\tau]$.
Let $D_d[0,\tau]$ be the Banach space of $d$-variate real valued cadlag functions endowed with a supremum norm $\pl\cdot\pl_{\infty}$ \citep{Neuhaus71}.
Let $L_1:Q({\cal M})\rightarrow D_d[0,\tau]$ and $L_2:G({\cal M})\rightarrow D_d[0,\tau]$. 
We assume that these loss functions and the canonical gradient map into functions in $D_d[0,\tau]$ with a sectional variation norm bounded by some universal finite constant (we will define sectional variation norm $\| .\|_v^*$ momentarily)
\begin{eqnarray}
M_1\equiv \sup_{P\in {\cal M}}\pl L_1(Q(P))\pl_v^* &<&\infty\nonumber, \\
 M_2\equiv \sup_{P\in {\cal M}}\pl L_2(G(P))\pl_v^*&<& \infty\nonumber, \\
 M_3\equiv \sup_{P\in {\cal M}}\pl D^*(P)\pl_v^*&<& \infty.
  \label{sectionalvarbound}
\end{eqnarray}
\begin{egbox}[(Treatment-specific mean)]
Under the previous stated assumptions,  the sectional variation norm of $W\rightarrow D^*(Q,G)(W,a,y)$ (for each $(a,y)\in \{0,1\}^2$) can be bounded in terms of the sectional variation norm of $W\rightarrow \bar{Q}(1,W)$ and $G$. Similarly, this same statement applies for $L(\bar{Q})$ and $L_2(G)$. As a consequence, the universal bounds (\ref{sectionalvarbound}) are finite.
\end{egbox}
 For a given function $F\in D_d[0,\tau]$, we define the sectional variation norm as follows. For a given subset $s\subset\{1,\ldots,d\}$,  let $F_s(x_s)=F(x_s,0_{-s})$ be the $s$-specific section of $F$ that sets the coordinates outside the subset $s$ equal to 0, where we used the notation $(x_s,0_{-s})$ for the vector whose $j$-th component equals $x_j$ if $j\in s$ and $0$ otherwise. 
 The sectional variation norm is now defined by
 \[
 \pl F\pl_v^*=| F(0)| +\sum_{s\subset\{1,\ldots,d\}}\int_{(0_s,\tau_s]}| dF_s(u_s)| ,\]
 where the sum is over all subsets $s$ of $\{1,\ldots,d\}$.
 Note that $\int_{(0_s,\tau_s]}| dF_s(u_s)| $ is the standard variation norm of the measure $dF_s$ generated by its $s$-specific section $F_s$ on the $| s|$-dimensional edge $(0_s,\tau_s]\times \{0_{-s}\}$ of the $d$-dimensional cube $[0,\tau]$. Thus, the sectional variation norm of $F$ is the sum of the variation norms  of $F$ itself and of all its $s$-specific sections $F_s$, plus that of the offset $|F(0)|$.
 We also note that any function $F\in D_d[0,\tau]$ with finite sectional variation norm (i.e., $\pl F\pl_v^*<\infty$) can be represented as follows \citep{Gill&vanderLaan&Wellner95}:
 \begin{equation}\label{Frepresentation}
 F(x)=F(0)+\sum_{s\subset\{1,\ldots,d\}}\int_{(0_s,x_s]} dF_s(u_s) .\end{equation}
As utilized in \citep{vanderLaan15} to define the HAL-MLE,  since $\int_{(0_s,x_s]}dF_s(u_s)=\int I_{u_s\leq x_s}dF_s(u_s)$, this representation shows that $F$ can be written as an infinitesimal linear combination of
 tensor product (over $s$) indicator basis functions $x\rightarrow I_{u_s\leq x_s}$ indexed by  a cut-off $u_s$, across all subsets $s$, where the coefficients in front of the tensor product indicator basis functions are equal to the infinitesimal increments $dF_s(u_s)$ of $F_s$ at $u_s$. This proves that this class of functions can be represented as a ''convex'' hull of the class of indicators basis functions, which proves that it is a Donsker class \citep{vanderVaart&Wellner96}.
 
For discrete measures $F_s$ this integral becomes a {\em finite} linear combination of such $| s|$-way indicator basis functions (where $| s|$ denotes the size of the set $s$). One could think of this representation of $F$ as a saturated model of a function $F$ in terms of tensor products of univariate indicator basis functions, ranging from products over singletons to product over the full set $\{1,\ldots,d\}$.  
For a function $f\in D_d[0,\tau]$,  we also define the supremum norm $\pl f\pl_{\infty}=\sup_{x\in [0,\tau]} |f(x)|$.


\subsubsection{General class of models for which parameter spaces for $Q$ and $G$ are Cartesian products of sets of cadlag functions with bounds on sectional variation norm}
Although the above bounds $M_1,M_2,M_3$ are the only relevant bounds for the asymptotic performance of the HAL-MLE and HAL-TMLE, for practical formulation of a model ${\cal M}$ one might prefer to state the sectional variation norm restrictions on the parameters $Q$ and $G$ themselves instead of on $L_1(Q)$ and $L_2(G)$. (In our formal results we will refer to such a model ${\cal M}$ as having the extra structure (\ref{calFmodel}) defined below, but, this extra structure is not needed, just as we can work with a general Donsker class as mentioned above.)

For that purpose, a model may assume that $Q=(Q_1,\ldots,Q_{K_1})$ for variation independent parameters $Q_k$  that are themselves $m_{1k}$-dimensional cadlag functions on $[0,\tau_{1k}]\subset \openr^{m_{1k}}_{\geq 0}$ with sectional variation norm bounded by some upper-bound $C_{Qk}^u$ and lower bound $C_{Qk}^l$, $k=1,\ldots,K_1$, and similarly  for $G=(G_1,\ldots,G_{K_2})$ with sectional variation norm bounds $C_{Gk}^u$ and $C_{Gk}^l$, $k=1,\ldots,K_2$. We define two parameters $Q_1$ and $Q_2$ are variation independent if $\{(Q_1(P), Q_2(P)): P \in \mathcal{M}\} = \{Q_1(P): P \in \mathcal{M}\} \otimes \{Q_2(P): P \in \mathcal{M}\}$ (i.e. tensor product of the parameter spaces of $Q_1$ and $Q_2$).
Typically, such a model would not enforce a lower bound on the sectional variation norm so that we have $C_{Qk}^l=C_{Gk}^l=0$.
Let $C_Q^u=(C_{Qk}^u:k=1,\ldots,K_1)$; $C_Q^l=(C_{Qk}^l:k=1,\ldots,K_1)$; and $C_Q=(C_Q^l,C_Q^u)$, and similarly we define $C_G^u$, $C_G^l$ and $C_G=(C_G^l,C_G^u)$.
Specifically, for such a class of models let
\begin{eqnarray*}
{\cal F}_{Qk}\equiv Q_k({\cal M}),\\
{\cal F}_{Gk}\equiv G_k({\cal M}),
\end{eqnarray*}
denote the parameter spaces for $Q_k$ and $G_k$,
and assume that these parameter spaces ${\cal F}_{jk}$ are contained in the class ${\cal F}_{jk}^{np}$ of $m_{jk}$-variate cadlag functions with sectional variation norm bounded from above by $C_{jk}^u$ and from below by $C_{jk}^l$, $k=1,\ldots,K_j$, $j \in \{Q, G\}$.
These bounds $C_Q^u=(C_{Qk}^u:k)$ and $C_G^u=(C_{Gk}^u:k)$ will then imply bounds $M_1,M_2,M_3$.
For such a model  $L_1(Q)$ and $L_2(G)$ would  be defined as sums of loss functions: $L_1(Q)=\sum_{k=1}^{K_1}L_{1k}(Q_k)$
and $L_2(G)=\sum_{k=1}^{K_2}L_{2k}(G_k)$. We also define the vector losses ${\bf L}_1(Q)=(L_{1k}(Q_k):k=1,\ldots,K_1)$, ${\bf L}_2(G)=(L_{2k}(G_k):k=1,\ldots,K_2)$,
and corresponding vector dissimilarities ${\bf d}_{01}(Q,Q_0)=(d_{01,k}(Q_k,Q_{k0}):k=1,\ldots,K_1)$ and ${\bf d}_{02}(G,G_0)=(d_{02,k}(G_k,G_{k0}):k=1,\ldots,K_2)$.

For example, the parameter space ${\cal F}_{jk}$ of $Q_k$ ($j=Q$) or $G_k$ ($j=G$) may be defined as \begin{equation}\label{calFmodel}
{\cal F}_{jk,A_{jk}}^{np}\equiv \{F\in {\cal F}_{jk}^{np}: dF_s(u_s)=I_{(s,u_s)\in A_{jk}}dF_s(u_s), s\subset\{1,\ldots,m_{jk}\}\},
\end{equation}  for some set $A_{jk}$ of possible values for $(s,u_s)$, $k=1,\ldots,K_j$, $j \in \{Q, G\}$,
where one evaluates this restriction on  $F$ in terms of the  representation (\ref{Frepresentation}).
Note that we used short-hand notation $g(x)=I_{x\in A} g(x)$ for $g$ being zero for $x\not \in A$.
We will make the convention that if $A$ excludes $\{0\}$, then it corresponds with assuming $F(0)=0$.

The subset  ${\cal F}_{Qk,A_{Qk}}^{np}$ of cadlag functions ${\cal F}_{Qk}^{np}$ with sectional variation norm between $C_{Qk}^l$ and $C_{Qk}^u$ further  restricts the support of these functions to a set $A_{Qk}$.
For example, $A_{Qk}$ might set $dF_s=0 $ for subsets $s$ of size larger than $3$ for all values $u_s\in (0_s,\tau_s]$, in which case the model assumes that the nuisance parameter $Q_k$ can be represented as a sum over all subsets $s$ of size $1,2$ and $3$ of a  function of the variables indicated by $s$.

In order to allow modeling of monotonicity (e..g, nuisance parameter $Q_k$ is an actual cumulative distribution function), we also allow that this set restricts $dF_s(u_s)\geq 0$ for all $(s,u_s)\in A_{jk}$. We will denote the latter parameter space with 
\begin{equation}\label{calFmodelplus}
{\cal F}_{jk,A_{jk}}^{np,+}=\{F\in {\cal F}_{jk}^{np}: dF_s(u_s)=I_{(s,u_s)\in A_{jk}}dF_s(u_s), dF_s\geq 0, F(0)\geq 0,\forall s\}.
\end{equation}
For the parameter space (\ref{calFmodelplus}) of monotone functions we allow that the sectional variation norm is known by setting $C_{jk}^u=C_{jk}^l$ (e.g, for the class of cumulative distribution functions we would have
$C_{jk}^u=C_{jk}^l=1$),  while for the parameter space (\ref{calFmodel}) of cadlag functions with sectional variation norm between $C_{jk}^l$ and $C_{jk}^u$ we assume 
$C_{jk}^l<C_{jk}^u$.



For the analysis of our proposed  nonparametric bootstrap sampling distributions  we do not assume this extra model structure that ${\cal F}_{jk}={\cal F}_{jk,A_{jk}}^{np}$ or ${\cal F}_{jk}={\cal F}_{jk,A_{jk}}^{np,+}$ for some set $A_{jk}$, $k=1,\ldots,K_j$, $j\in\{Q, G\}$.
In the sequel we will refer to a model with this extra structure as a model satisfying (\ref{calFmodel}), even though we include the case (\ref{calFmodelplus}). 
All our formal results apply without this  extras model structure (and also for any other uniform Donsker class as mentioned above), but it just happens to represent a natural model structure for establishing the sectional variation norm bounds (\ref{sectionalvarbound}) on $L_1(Q)$, $L_2(G)$, and $D^*(Q,G)$, and for computing HAL-MLEs.
The key  practical benefit of this extra model structure is that the implementation of the HAL-MLE for such a parameter space ${\cal F}_{jk,A_{jk}}^{np}$ corresponds with fitting a linear combination of indicator basis functions of the form
$I_{u_s\leq x_s}$ (indexed by a subset $s$ and knot-point $u_s$)
under the sole constraint that the sum of the absolute value of the coefficients is bounded by $C_{jk}^l$ and $C_{jk}^u$, and possibly that the coefficients are non-negative, where the set $A_{jk}$ implies the set of indicator basis functions that are included.
Specifically, in the case that the nuisance parameter is a conditional mean or conditional probability we can compute the  HAL-MLE with standard lasso linear or logistic regression software \citep{Benkeser&vanderLaan16}. Therefore, this restriction on our set of models also allows straightforward computation of its HAL-MLEs,  corresponding HAL-TMLE, and their bootstrap analogues. 

A typical statistical model assuming the extra structure (\ref{calFmodel}) would be of the form ${\cal M}=\{P: Q_{k_1}(P)\in {\cal F}_{Qk_1,A_{Qk_1}}^{np},G_{k_2}(P)\in {\cal F}_{Gk_2,A_{Gk_2}}^{np},k_1,k_2\}$ indexed by the support sets $((A_{Qk_1},A_{Gk_2}):k_1,k_2)$ and the sectional variation norm  bounds $((C_{jk}^l,C_{jk}^u):j,k)$, but the model ${\cal M}$ might  include additional restrictions on $P$ as long as the parameter spaces of these nuisance parameters equal these sets ${\cal F}_{jk_j,A_{jk_j}}^{np}$ or ${\cal F}_{jk_j,A_{jk_j}}^{np,+}$.

\begin{remark}[Creating  parameter spaces of type (\ref{calFmodel}) or (\ref{calFmodelplus})]
In our first example we have a nuisance parameter $\bar{G}(W)=E_P(A\mid W)$ that is not just assumed to be cadlag and have bounded sectional variation norm but is also bounded between $\delta$ and $1-\delta$ for some $\delta>0$. This means that the parameter space for this ${G}$ is not exactly of type (\ref{calFmodel}). 
This is easily resolved by, for example, reparameterizing $\bar{G}(W)=\mbox{expit}(G(W))$ where $G$ can be any cadlag function with sectional variation norm bounded by some constant $C^u$.
The bound $C^u$ implies automatically a supremum norm bound on $G$, and thereby that $\delta<\bar{G}<1-\delta$ for some $\delta=\delta(C^u)>0$. One now defines the nuisance parameter as $G$. Similarly, such a parametrization can be applied to $E(Y\mid A,W)$ and to the density in our second example.
 These just represent a few examples showcasing that one can reparametrize the natural nuisance parameters in terms of nuisance parameters that have a parameter space of the form (\ref{calFmodel}) or (\ref{calFmodelplus}). These representations are actually natural steps for the implementation of the HAL-MLE since they allow us now to minimize the empirical risk over a generalized linear model with the sole constraint that the sum of absolute value of coefficients is bounded  (and possibly coefficients are non-negative). 
\end{remark}

\subsubsection{Bounding the exact second-order remainder in terms of loss-based dissimilarities}
Let \[R_2(P,P_0)=R_{20}(Q,G,Q_0,G_0)\]
for some mapping $R_{20}()=R_{2P_0}()$ possibly indexed by $P_0$. 
We assume the following upper bound:
\begin{equation}\label{boundingR2}| R_2(P,P_0)| =| R_{20}(Q,G,Q_0,G_0)| \leq f({\bf d}_{01}^{1/2}(Q,Q_0),{\bf d}_{02}^{1/2}(G,G_0))\end{equation}
for some function $f:\openr^K_{\geq 0}\rightarrow\openr_{\geq 0}$, $K=K_1+K_2$, of the form $f(x)=\sum_{i,j} a_{ij} x_ix_j$, a quadratic polynomial with positive coefficients $a_{ij}\geq 0$.  In all our examples, one simply uses the Cauchy-Schwarz inequality to bound $R_{20}(P,P_0)$ in terms of $L^2(P_0)$-norms of $Q_{k_1}-Q_{k_10}$ and $G_{k_2}-G_{k_20}$, and subsequently one relates these  $L^2(P_0)$-norms to its loss-based dissimilarities $d_{01,k_1}(Q_{k_1},Q_{k_10})$ and $d_{02,k_2}(G_{k_2},G_{k_20})$, respectively. This bounding step will also rely on an assumption  that denominators in $R_{20}(P,P_0)$ are uniformly bounded away from zero. This type of assumption that guarantees uniform bounds on $D^*(Q,G)$ and on $R_{20}(Q,G,Q_0,G_0)$ is often referred to as a strong positivity assumption since it requires that the data density has a certain type of support relevant for the target parameter $\Psi$, and that the data density is uniformly bounded away from zero on that support. 
In the treatment specific mean example, a common case where the strong positivity assumption does not hold is if $G_0(A=1|W) = 0$ for a some value of $W$.


\subsubsection{Continuity of efficient influence curve as function of $P$ at $P_0$}
We also assume that if the rates of convergence of $d_{01}(Q_n,Q_0)$ and $d_{02}(G_n,G_0)$ translate in the same rate of convergence of $P_0\{D^*(Q_n,G_n)-D^*(Q_0,G_0)\}^2$. This is guaranteed by the following upper bound:
\begin{equation}\label{contDstar}
P_0\{D^*(Q,G)-D^*(Q_0,G_0)\}^2 \leq f({\bf d}_{01}^{1/2}(Q,Q_0),{\bf d}_{02}^{1/2}(G,G_0))\end{equation}
for some function $f:\openr^K_{\geq 0}\rightarrow\openr_{\geq 0}$, $K=K_1+K_2$, of the form $f(x)=\sum_{i,j} a_{ij} x_ix_j$, a quadratic polynomial with positive coefficients $a_{ij}\geq 0$. 

\subsection{HAL-MLEs of nuisance parameters}
\sloppy We estimate $Q_0,G_0$ with HAL-MLEs $Q_n,G_n$ satisfying (with probability tending to 1)
\begin{eqnarray*}
P_n L_1(Q_n)&\leq &P_n L_1(Q_0),\\
P_n L_2(G_n)&\leq &P_n L_2(G_0).
\end{eqnarray*}
For example, $Q_n$ might be defined as the actual minimizer $Q_n=\argmin_{Q\in Q({\cal M})}P_n L_1(Q)$.
If $Q$ has multiple components and the loss function is a corresponding sum loss function, then  these HAL-MLEs correspond with separate HAL-MLEs for each component. 
We have the following previously established result from Lemma 3 in \cite{vanderLaan15} for these HAL-MLEs. We represent estimators as mappings on the nonparametric model ${\cal M}^{np}$ containing all possible realizations of the empirical measure $P_n$.
\begin{lemma}\label{lemma2}
\textbf{(Lemma 3 from \cite{vanderLaan15})}
Let $O\sim P_0\in {\cal M}$. Let $Q:{\cal M}\rightarrow Q({\cal M})$ be a function valued parameter and let $L:Q({\cal M})\rightarrow D_d[0,\tau]$ be a loss function so that $Q_0\equiv Q(P_0)=\argmin_{Q\in Q({\cal M})}P_0L(Q)$. 
 Let $\hat{Q}:{\cal M}^{np}\rightarrow Q({\cal M})$ define an estimator  $Q_n\equiv \hat{Q}(P_n)$ so that $P_n L_1(Q_n)=\min_{Q\in Q({\cal M})}P_n L(Q)$ or $P_n L_1(Q_n)\leq P_n L_1(Q_0)$. 
 Let $d_0(Q,Q_0)=P_0L(Q)-P_0L(Q_0)$ be the loss-based dissimilarity. 
 Then,
 \[
 d_0(Q_n,Q_0)\leq -(P_n-P_0)\{L(Q_n)-L(Q_0)\}. 
 \]
 If $\sup_{Q\in Q({\cal M})}\pl L(Q)\pl_v^*<\infty$, and (\ref{M20bounds}) holds for $L_1(Q)$, then 
 \[
E_0 d_0(Q_n,Q_0)=O(n^{-1/2-\alpha/4}),\]
where $\alpha$ is defined as in  (\ref{alphad}) for class $\{L_1(Q):Q\in Q({\cal M})\}$.
 \end{lemma}
 Application of this general lemma proves that $d_{01}(Q_n,Q_0)=O_P(n^{-1/2-\alpha/4})$ and $d_{02}(G_n,G_0)=O_P(n^{-1/2-\alpha/4})$. 

One can add restrictions to the parameter space $Q({\cal M})$ over which one minimizes in the definition of $Q_n$ and $G_n$ as long as one guarantees that, with probability tending to 1, $P_nL_1(Q_n)\leq P_nL_1(Q_0)$ and $P_n L_2(G_n)\leq P_n L_2(G_0)$. For example, in a model ${\cal M}$ with extra structure (\ref{calFmodel}) this allows one to use a data dependent upper bound $C_{Qn}^u\leq C^u_1$ on the sectional variation norm  in the definition of $Q_n$ if we know that $C_{Qn}^u$ will be larger than the true $C_{Q0}^u=\pl Q_0\pl_v^*$ with probability tending to 1.

\subsection{HAL-TMLE}
Consider a finite dimensional local  least favorable model $\{Q_{n,\epsilon}:\epsilon\}\subset Q({\cal M})$ through $Q_n$ at $\epsilon =0$ so that  the linear span of the components of $\frac{d}{d\epsilon}L_1(Q_{n,\epsilon})$ at $\epsilon =0$ includes $D^*(Q_n,G_n)$. 
Let $Q_n^*=Q_{n,\epsilon_n}$ for $\epsilon_n=\argmin_{\epsilon}P_n L_1(Q_{n,\epsilon})$. We assume that this one-step TMLE  $Q_n^*$ already satisfies
\begin{equation}\label{efficeqn}
r_n\equiv P_n D^*(Q_n^*,G_n)=o_P(n^{-1/2}).\end{equation}
Since $d_{01}(Q_n,Q_0)=o_P(n^{-1/2})$ we will  have that $\epsilon_n=o_P(n^{-1/4})$, and $\epsilon_n$ solves its score equation $\frac{d}{d\epsilon_n}P_n L_1(Q_{n,\epsilon_n})=0$, which, in first order, equals its score equation $P_n D^*(Q_{n,\epsilon_n},G_n)$ at $\epsilon=0$ (with a second order remainder $O(\epsilon_n^2)=o_P(n^{-1/2})$). 
This basic argument allows one to prove that (\ref{efficeqn}) holds under the assumption $d_{01}(Q_n,Q_0)=o_P(n^{-1/2})$ and regularity conditions, as formally shown in the Appendix of \citep{vanderLaan15}.
 Alternatively,  one could use the one-dimensional canonical universal least favorable model satisfying $\frac{d}{d\epsilon}L_1(Q_{n,\epsilon})=D^*(Q_{n,\epsilon},G_n)$ at each $\epsilon$ (see our second example in Section \ref{sectfive}). In that case, the efficient influence curve equation (\ref{efficeqn}) is solved exactly with the one-step TMLE: i.e.,  $r_n=0$ \citep{vanderLaan&Gruber15}. 
 The HAL-TMLE of $\psi_0$ is the plug-in estimator $\psi_n^*=\Psi(Q_n^*)$. 
 In the context of model structure (\ref{calFmodel}) (or (\ref{calFmodelplus})), we will also refer to   this estimator as the HAL-TMLE$(C^u)$ to indicate its dependence on the specification of the bounds $C^u=(C_Q^u,C_G^u)$ on the sectional variation norms of the components of $Q$ and $G$
 
Lemma \ref{lemmathalmle} in Appendix A proves that $d_{01}(Q_{n,\epsilon_n},Q_0)$ converges at the same rate as $d_{01}(Q_n,Q_0)=O_P(n^{-1/2-\alpha/4})$ (see (\ref{thalmle})). This also implies this result for any $K$-th step TMLE with $K$ fixed.
The advantage of  a one-step or $K$-th step TMLE is that it is always well defined, and it easily follows that it converges at the same rate as the initial $Q_n$ to $Q_0$. In addition, for these closed form TMLEs it is also guaranteed that the sectional variation norm of $Q_n^*$ remains universally  bounded. The latter is important for the Donsker class condition for asymptotic efficiency of the HAL-TMLE, but the Donsker class condition could be avoided by using a   cross-validated HAL-TMLE that relies on sample splitting \citep{vanderLaan&Rose11}.

Assuming extra model structure (\ref{calFmodel}), since we apply the least favorable submodel to an HAL-MLE $Q_n$ that is likely having the maximal allowed $C^u_1$ sectional variation norm, the following remark is in order. We suggest to  simply extend the statistical model by enlarging the sectional variation norm bounds to $C^u_1+\delta$ for some $\delta>0$, even though the original bounds $C^u_1$ are still used in  the definition of the HAL-MLEs. This increase in statistical model {\em does not change} the canonical gradient at $P_0$ (known to be an element of the interior of original model), while now a  least favorable submodel through the HAL-MLE is allowed to enlarge the sectional variation norm. This makes the construction of a least favorable submodel easier by not having to worry to constrain the sectional variation norm. Since the HAL-MLE $Q_n$ has the maximal allowed uniform sectional variation norm $C_Q^u$, and $Q_n$ is consistent, the sectional variation norm of  the TMLE $Q_n^*=Q_{n,\epsilon_n}$ will now be slightly larger, and asymptotically approximate $C^u_1$. 
Either way, with the slightly enlarged definition of ${\cal M}$, we have $\{Q_{n,\epsilon}:\epsilon\}\subset{\cal M}$ so that the assumption (\ref{sectionalvarbound}) guarantees that $\pl L_1(Q_{n,\epsilon_n})\pl_v^*$ is bounded by a universal constant. 

\begin{egbox}[(Treatment-specific mean)]
Condition (\ref{boundingR2}) holds by applying the Cauchy-Schwarz inequality, and using $G>\delta>0$ for some $\delta>0$.    The HAL-MLEs $\bar{Q}_n$ and $G_n$ of $\bar{Q}$ and $G$, respectively, can be computed with a lasso-logistic regression estimator with large (approximately $n 2^d$) number of indicator basis functions (see our example section for more details), where we can select the $L^1$-norm of the coefficient vector with cross-validation. The least favorable submodel through $\bar{Q}_n$ is given by
    \begin{align}
    \mbox{logit} \bar{Q}_{n,\varepsilon} = \mbox{logit}\bar{Q}_n + \varepsilon C(G_n),
    \end{align}
where $C(G_n)(A,W) \triangleq A/G_n(W)$. Let $\varepsilon_n \triangleq \argmin_\varepsilon P_n L_{1}(Q_{n,\varepsilon})$, which is thus computed with a simple univariate logistic regression MLE, using as off-set $\mbox{logit}\bar{Q}_n$. This defines the TMLE $\bar{Q}_n^*=\bar{Q}_{n,\epsilon_n}$. Recall that $Q_{W,n}$ is already an NPMLE so that a TMLE-update based on a log-likelihood loss and local least favorable submodel (i.e., with score $\bar{Q}_n(W)-\Psi(Q_n)$, will not change this estimator. Let $Q^*_n=(Q_{W,n},\bar{Q}_n^*)$.
The HAL-TMLE  of $\psi_0$ is the plug-in estimator $\psi^*_n \triangleq \Psi(Q_n^*) =\frac{1}{n}\sum_{i=1}^n \bar{Q}_n^*(1,W_i)$.
\end{egbox}

\subsection{Asymptotic efficiency theorem for HAL-TMLE and CV-HAL-TMLE}
Lemma \ref{lemma2} establishes that $d_{01}(Q_n,Q_0)$ and $d_{02}(G_n,G_0)$ are $O_P(n^{-1/2-\alpha/4})$. Lemma \ref{lemmathalmle} in Appendix A proves that also ${ d}_{01}(Q_n^*,Q_0)=O_P(n^{-1/2-\alpha/4})$. Combined with (\ref{boundingR2}), this shows that the second-order term $R_{20}(Q_n^*,G_n,Q_0,G_0)=O_P(n^{-1/2-\alpha/4})$.

We have the following identity for the HAL-TMLE:
\begin{eqnarray}
\Psi(Q_n^*)-\Psi(Q_0)&=&(P_n-P_0)D^*(Q_n^*,G_n)+R_{20}(Q_n^*,G_n,Q_0,G_0)+r_n \label{exactexpansiontmle_first_equality}\\
&=&(P_n-P_0)D^*(Q_0,G_0)+(P_n-P_0)\{D(Q_n^*,G_n)-D^*(Q_0,G_0)\}\nonumber \\
&&+R_{20}(Q_n^*,G_n,Q_0,G_0)+r_n.\label{exactexpansiontmle}
\end{eqnarray}
The second term on the right-hand side is $O_P(n^{-1/2-\alpha/4})$ following the same empirical process theory proof as  Theorem 1 in \cite{vanderLaan17} using the continuity condition (\ref{contDstar}) on $D^*$.
Thus, this proves the following asymptotic efficiency theorem.

\begin{theorem}\label{thefftmle}
Consider the statistical model ${\cal M}$ and target parameter $\Psi:{\cal M}\rightarrow\openr$ s
satisfying
(\ref{M20bounds}), (\ref{sectionalvarbound}), (\ref{boundingR2}), (\ref{contDstar}).
Let $Q_n,G_n$ be the above defined HAL-MLEs, where $d_{01}(Q_n,Q_0)$ and $d_{02}(G_n,G_0)$ are $O_P(n^{-1/2-\alpha/4})$. Let $Q_n^*=Q_{n,\epsilon_n}$ be the one-step TMLE-update  according to a submodel $\{Q_{n,\epsilon}:\epsilon\}\subset {\cal M}$ solving the efficient influence curve equation such that (\ref{efficeqn}) holds.

Then the HAL-TMLE $\Psi(Q_n^*)$ of $\psi_0$ is asymptotically efficient:
\begin{equation}\label{efficiencyhaltmle}
\Psi(Q_n^*)-\Psi(Q_0)=P_nD^*(Q_0,G_0)+O_P(n^{-1/2-\alpha/4}).\end{equation}
\end{theorem}
We remind the reader that the condition (\ref{sectionalvarbound}), stating that the loss functions and canonical gradient are contained in class of cadlag functions with a universal bound on the sectional variation norm, can be replaced by a general Donsker class condition (\ref{alphad}). We also remark that this Theorem \ref{thefftmle}  trivially generalizes to any rate of convergence for $d_{01}(Q_n,Q_0)$ and $d_{02}(G_n,G_0)$ by simply setting the remainder term in (\ref{efficiencyhaltmle}) equal to this same rate. Due to a recent new result in \citep{Bibaut&vanderLaan19} for the HAL-MLE, utilizing an improved recently published covering number bound,  under specified conditions, we now even have $d_{01}(Q_n,Q_0)=O_P(n^{-2/3}(\log n)^d)$ and $d_{02}(G_n,G_0)=O_P(n^{-2/3}(\log n)^d)$.  So, applying this result would yield (\ref{efficiencyhaltmle}) with remainder $O_P(n^{-2/3}(\log n)^d)$.

\subsubsection{Wald type confidence interval}
A first order asymptotic 0.95-level confidence interval is given by $\psi_n^*\pm 1.96 \sigma_n/n^{1/2}$ where $\sigma_n^2=P_n \{D^*(Q_n^*,G_n)\}^2$ is a consistent estimator of $\sigma^2_0=P_0\{D^*(Q_0,G_0)\}^2$.
Clearly, this first order confidence interval ignores the exact remainder $\tilde{R}_{2n}$ in the exact expansion $\Psi(Q_n^*)-\Psi(Q_0)=(P_n-P_0)D^*(Q_0,G_0)+\tilde{R}_{2n}$ as presented in (\ref{exactexpansiontmle}):
\begin{equation}\label{exactremainder_tmle}
\tilde{R}_{2n}\equiv R_{20}(Q_n^*,G_n,Q_0,G_0)+(P_n-P_0)\{D^*(Q_n^*,G_n)-D^*(Q_0,G_0)\}+r_n.\end{equation}
 
Let's consider the extra model structure (\ref{calFmodel}). The asymptotic efficiency proof above of the HAL-TMLE$(C^u)$ relies on the HAL-MLEs $(Q_{n,C_Q^u},G_{n,C_G^u})$ converging to the true $(Q_0,G_0)$ at rate faster than $n^{-1/4}$, and their sectional variation norm being uniformly bounded from above by $C^u=(C_Q^u,C_G^u)$. Both of these conditions are still known to hold for the CV-HAL-MLE $(Q_{n,C_{Qn}},G_{n,C_{Gn}})$ in which the constants $(C_Q,C_G)$ are selected with the cross-validation selector $C_n=(C_{Qn},C_{Gn})$ \citep{vanderLaan15}. This follows since the cross-validation selector is asymptotically equivalent to the oracle selector, thereby guaranteeing that $C_n$ will exceed the sectional variation norm of the true $(Q_0,G_0)$ with probability tending to 1. Typically, one will only data adaptively select $C^u$, while keeping $C^l=(C_Q^l,C_G^l)$ at its known lower bound.
Therefore, we  have that this CV-HAL-TMLE is also asymptotically efficient. Of course, this CV-HAL-TMLE is more practical and powerful than the HAL-TMLE at an apriori specified $C=(C_Q,C_G)=(C_Q^u,C_Q^l,C_G^u,C_G^l)$  since it adapts the choice of bounds $C=(C_Q,C_G)$ to  the true sectional variation norms $C_0=(C_{Q0},C_{G0})$ for $(Q_0,G_0)$.
 
 For simplicity, in the next theorem we focus on data adaptive selection of $C^u$ only.
 \begin{theorem}\label{thefftmlecv}
 Consider the setting of Theorem \ref{thefftmle}, but with the extra model structure (\ref{calFmodel}). 
 Let $C_{Q0}^u=\pl Q_0\pl_v^*$, $C_{G0}^u=\pl G_0\pl_v^*$.
 Suppose that $C_Q^u$ and $C_G^u$  that define the HAL-MLEs $Q_n=Q_{n,C_Q^u}$ and $G_n=G_{n,C_G^u}$ are replaced by data adaptive selectors $C_{Qn}^u$ and $C_{Gn}^u$ for which 
 \begin{equation}
 P_0(C_{Q0}^u\leq C_{Qn}^u\leq C_Q^u,C_{G0}^u\leq C_{Gn}^u\leq C_G^u)\rightarrow 1,\mbox{ as $n\rightarrow\infty$.}\label{Cn}
 \end{equation}
  Then, under the same assumptions as in Theorem \ref{thefftmle}, the TMLE
 $\Psi(Q_n^*)$,  using $Q_n=Q_{n,C_{Qn}^u}$ and $G_n=G_{n,C_{Gn}^u}$ as initial estimators,  is asymptotically efficient. 
 \end{theorem}
 In general, when the model ${\cal M}={\cal M}(C)$ is defined by global  constraints $C$, then one should use cross-validation to select these constraints $C$, which will only improve the performance of the initial estimators and corresponding TMLE, due to its asymptotic equivalence with the oracle selector. So our model ${\cal M}$  satisfying (\ref{sectionalvarbound}) and the extra structure (\ref{calFmodel}) might have more global constraints beyond $C^u=(C_Q^u,C_G^u)$ and these could then also be selected with cross-validation resulting in a CV-HAL-MLE and corresponding HAL-TMLE (see also our two examples).

\section{The nonparametric bootstrap for the HAL-TMLE}\label{sectthree}
Let $O_1^{\#},\ldots,O_n^{\#}$ be $n$ i.i.d. draws from the empirical measure $P_n$.
Let $P_n^{\#}$ be the empirical measure of this bootstrap sample. 

\subsection{Definition of bootstrapped HAL-MLEs for model with extra structure (\ref{calFmodel})}
In this subsection, we will assume the extra structure (\ref{calFmodel}) so that our parameter spaces for $Q$ and $G$ consists of cadlag functions with a universal bound $C^u$ on the sectional variation norm, thereby allowing us specific computational friendly definitions of the bootstrapped HAL-MLEs.  We generalize the definition of $Q$ being absolutely continuous w.r.t. $Q_n$:
$Q\ll Q_n$.
\begin{definition}\label{defabscont}
Recall the representation (\ref{Frepresentation}) for a multivariate real valued cadlag function $F$ in terms of its sections $F_s$.
Assume the extra model structure (\ref{calFmodel}) on ${\cal M}$.
We will say that $Q_k$ is absolutely continuous w.r.t. $Q_{k,n}$ if for each subset $s\subset\{1,\ldots,m_{1k}\}$, its $s$-specific section $Q_{k,s}$ defined by
$u_s\rightarrow Q_{k}(u_s,0_{-s})$ is absolutely continuous w.r.t. $Q_{n,k,s}$ defined by $u_s\rightarrow Q_{n,k}(u_s,0_{-s})$.
 We use the notation $Q_k\ll Q_{n,k}$. In addition, we use the notation $Q\ll Q_n$ if $Q_k\ll Q_{n,k}$ for each component $k\in \{1,\ldots,K_1\}$.
 Similarly, we use this notation $G\ll G_n$ if $G_k\ll G_{n,k}$ for each component $k\in \{1,\ldots,K_2\}$.
 \end{definition}
 In practice, the HAL-MLE $Q_n=\argmin_{Q\in Q({\cal M})}P_n {\bf L}_1(Q)$ is attained (or simply defined as a minimum among all discrete measures with fine enough selected support) by a discrete measure  $Q_n$ so that it can be computed by  minimizing the empirical risk over a large linear combination of indicator basis functions (e.g., $2^{m_{1k}} n$  for $Q_{nk}$) under the constraint that the sum of the absolute value of the coefficients is bounded by the specified constant $C_Q$ \citep{Benkeser&vanderLaan16}.
In that case, the constraint $Q\ll Q_n$ states that $Q$ is a linear combination of the indicator basis functions that had a non-zero coefficient in $Q_n$.

 Let  \begin{eqnarray*}
 Q_n^{\#}&=&\argmin_{Q\in Q({\cal M}),Q\ll Q_n,\pl Q\pl_v^*\leq \pl Q_n\pl_v^*}P_n^{\#}L_1(Q),\\ 
 G_n^{\#}&=&\argmin_{G\in G({\cal M}),G\ll G_n,\pl G\pl_v^*\leq \pl G_n\pl_v^*}P_n^{\#}L_2(G)
 \end{eqnarray*}
 be the corresponding HAL-MLEs of $Q_n=\argmin_{Q\in Q({\cal M})}P_nL_1(Q)$ and $G_n=\argmin_{G\in G({\cal M})}P_nL_2(G)$ based on the bootstrap sample. Here $P_n^{\#}$ is the empirical probability measure that puts mass $1/n$ on each observation $O_i^{\#}$ from a bootstrap sample $O_1^{\#},\ldots,Q_n^{\#}$ represnting $n$ i.i.d. draws from $P_n$.
Since our results on the rates of convergence of $Q_n^{\#}$ and $G_n^{\#}$ to $Q_n$ and $G_n$ only rely on $P_n^{\#}L_1(Q_n^{\#})\leq P_n^{\#}L_1(Q_n)$ (and similarly for $G_n^{\#})$, the additional restrictions $Q\ll Q_n$ and $\pl Q\pl_v^*\leq \pl Q_n\pl_v^*$ are appropriate theoretically.
In addition, the extra restriction $Q\ll Q_n$ makes the computation of the HAL-MLE on the bootstrap sample much faster than the HAL-MLE $Q_n$ based on the original sample, so that enforcing this extra constraint is only beneficial from a computational point of view. That is, the computation of $Q_n^{\#}$ only involves minimizing the empirical risk w.r.t. $P_n^{\#}$ over the coefficients that were  non-zero in the $Q_n$-fit. Given our experience that a typical HAL-MLE fit has around $n$ non-zero coefficients, this makes the calculation of $Q_n^{\#}$ across many bootstrap samples computationally feasible. Additionally in Appendix \ref{AppendixF}, we include two empirical simulation results on the number of non-zero coefficients as a function of sample size.

The above bootstrap distribution depends on the bounds $C=(C_Q,C_G)$ enforced in the HAL-MLEs $(Q_n,G_n)$. One possible choice is to set $C=(C_Q,C_G)$ equal to the cross-validation selector $C_{n,cv}$, where we typically only adaptively select the upper bound $C^u$ so that $C_{n,cv}=(C_{n,cv}^u,C^l)$.  In the next section we discuss an alternative (so called plateau) selector $C_{n}^u$ for  $C^u$ that aims to improve finite sample coverage. Either way, in the bootstrap distribution the choice $C$ ($=C_{n,cv}$ or $=C_n$) is treated as fixed, although we will evaluate the bootstrap distribution for a range of $C$ values to determine the plateau selector $C_n^u$.

\subsection{Definition of bootstrapped HAL-MLE in general.}
In general, $Q_n^{\#}\in Q({\cal M})$ and $G_n^{\#}\in G({\cal M})$ have to be defined as estimators of $Q_n$ and $G_n$ based on the bootstrap sample $P_n^{\#}$ satisfying that, with probability tending to 1 (conditional on $P_n$), $P_n^{\#}L_1(Q_n^{\#})\leq P_n^{\#} L_1(Q_n)$ and $P_n^{\#}L_2(G_n^{\#})\leq P_n^{\#}L_2(G_n)$. 
For example,  $Q_n^{\#}=\argmin_{Q\in Q({\cal M})}P_n^{\#}L_1(Q)$
and $G_n^{\#}=\argmin_{G\in G({\cal M})}P_n^{\#}L_2(G)$, but one is  allowed to add restrictions to the parameter space over which one minimizes $Q\rightarrow P_n^{\#}L_1(Q)$ as long as this space still includes $Q_n$ with probability tending to 1 (and similarly for $G_n^{\#}$).

\subsection{Bootstrapped HAL-TMLEs}
Let $\epsilon_n^{\#}=\argmin_{\epsilon}P_n^{\#}L_1(Q_{n,\epsilon}^{\#})$ be the one-step TMLE update of $Q_n^{\#}$ based on the least favorable submodel $\{Q_{n,\epsilon}^{\#}:\epsilon\}$ through $Q_n^{\#}$ at $\epsilon =0$ with score $D^*(Q_n^{\#},G_n^{\#})$ at $\epsilon =0$.  Let $Q_n^{\#*}=Q_{n,\epsilon_n^{\#}}^{\#}$ be the TMLE update which is assumed to solve 
\begin{equation}\label{efficeqnb}
r_n^{\#}\equiv |P_n^{\#}D^*(Q_n^{\#*},G_n^{\#})| =o_{P_n}(n^{-1/2}),
\end{equation}
 conditional on $(P_n:n\geq 1)$ (just like $r_n=o_P(n^{-1/2})$). Let $\Psi(Q_n^{\#*})$ be the resulting TMLE of $\Psi(Q_n^*)$ based on this nonparametric bootstrap sample.  Let $\sigma_n^2$ be an estimate of the asymptotic variance $\sigma^2_0=P_0 D^*(Q_0,G_0)^2$, such as $\sigma_n^2=P_n D^*(Q_n,G_n)^2$. Let $\sigma_n^{\#2}$ be this estimator applied to $P_n^{\#}$. 
We estimate the finite sample distribution of $n^{1/2}(\Psi(Q_n^*)-\Psi(Q_0))/\sigma_n$ with the sampling distribution of $Z_n^{1,\#}\equiv n^{1/2}(\Psi(Q_n^{\#*})-\Psi(Q_n^*))/\sigma_n^{\#}$, conditional on $P_n$.
Let $\Phi_n^{\#}(x)=P(n^{1/2}(\Psi(Q_n^{\#*})-\Psi(Q_n^*))/\sigma_n^{\#}\leq x\mid P_n)$ be the cumulative distribution of this bootstrap sampling distribution.
So a bootstrap based 0.95-level confidence interval for $\psi_0$ is given by \[
[\psi_n^{*}+q_{0.025,n}^{\#}\sigma_n/n^{1/2},\psi_n^*+q_{0.975,n}^{\#}\sigma_n/n^{1/2} ],\]
 where 
$q_{p,n}^{\#}=\Phi_n^{\#-1}(p)$ is the $p$-th quantile of this bootstrap distribution. 
We note that the upper bounds $\pl Q_n\pl_v^*$ and $\pl G_n\pl_v^*$ on the sectional variation norms of $Q_n^{\#}$ and $G_n^{\#}$, or equivalently, the upper bounds $C_{Qn}^u$ and $C_{Gn}^u$ in the definition of the HAL-MLEs $Q_n$ and $G_n$, will impact the values of these quantiles
$q_{0.025,n}^{\#}$ and $q_{0.975,n}^{\#}$. That is, the larger these values, the larger the finite dimensional models for $Q_n$ and $G_n$ implied by their non-zero coefficients, and thereby the larger the variation of the resulting TMLE $\Psi(Q_n^{\#*})$. Our results apply for any data adaptive selector $C_n$ satisfying that, with probability tending to 1, $C_{Qn}^u$ is larger than $\pl Q_0\pl_v^*$ and smaller than $C_Q^u$, and similarly for $C_{Qn}^u$. However, clearly, the finite sample coverage of the resulting bootstrap confidence interval is affected by the precise choice $C_n^u=(C_{Qn}^u,C_{Gn}^u)$.


We now want to prove that $n^{1/2}(\Psi(Q_n^{\#*})-\Psi(Q_n^*))$, conditional on $P_n$, converges in distribution to $N(0,\sigma^2_0)$, and thereby also that $\Phi_n^{\#}$ converges to  the cumulative distribution function of limit distribution $N(0,1)$. Importantly, this nonparametric bootstrap confidence interval could potentially dramatically improve the coverage relative to using the first order Wald-type confidence interval since this bootstrap distribution is estimating the variability of the full-expansion of the TMLE, including the exact remainder $\tilde{R}_{2n}$.

In the next subsection we show that the nonparametric bootstrap works for the HAL-MLEs $Q_n$ and $G_n$. Subsequently, not surprisingly,  we can show that this also establishes that the bootstrap works for the one-step TMLE $Q_n^*$ ($K$-th step TMLE for fixed $K$). This provides then the basis for proving that the nonparametric bootstrap is consistent for the HAL-TMLE.

\subsection{Nonparametric bootstrap for HAL-MLE}

The following theorem establishes that the bootstrap HAL-MLE $Q_n^{\#}$  estimates $Q_n$ as well, w.r.t.\ an empirical loss-based dissimilarity $d_{n1}(Q_n^{\#},Q_n)=P_n L_1(Q_n^{\#})-P_nL_1(Q_n)$, as
$Q_n$  estimates  $Q_0$  with respect to $d_{01}(Q_n,Q_0)=P_0L_1(Q_n)-P_0L_1(Q_0)$. 
In fact, we even have $d_{01}(Q_n^{\#},Q_0)=O_P(n^{-1/2-\alpha/4})$.
The analogue results apply to $G_n^{\#}$. 

\begin{theorem}\label{thnpbootmle} 
Assume (\ref{M20bounds}) and (\ref{sectionalvarbound}).

{\bf Definitions:}
Let $d_{n1}(Q,Q_n)=P_n \{L_1(Q)-L_1(Q_n)\}$ be the loss-based dissimilarity at the empirical measure, where $Q_n$ is an HAL-MLE of $Q_0$ satisfying $P_n L_1(Q_n)\leq P_n L_1(Q_0)$. 
Similarly, let $d_{n2}(G,G_n)=P_n \{L_2(G)-L_2(G_n)\}$ be the loss-based dissimilarity at the empirical measure, where $G_n$ is an HAL-MLE of $G_0$ satisfying $P_n L_2(G_n)\leq P_n L_2(G_0)$. 

{\bf Conclusion:}
Then,
\[
d_{n1}(Q_n^{\#},Q_n)=O_P(n^{-1/2-\alpha/4})\mbox{ and } d_{n2}(G_n^{\#},G_n)=O_P(n^{-1/2-\alpha/4}).\]
We also have 
 \[ d_{01}(Q_n^{\#},Q_0)=O_P(n^{-1/2-\alpha/4})\mbox{ and } d_{02}(G_n^{\#},G_0)=O_P(n^{-1/2-\alpha/4}).\]
 {\bf Bootstrapping HAL-MLE$(C)$ at $C^u=C_n^u$ for model with extra structure (\ref{calFmodel}):}
 This result also applies to the case that $C^u=(C_Q^u,C_G^u)$ in definition of HAL-MLEs $(Q_n,G_n)$  is replaced by a data adaptive choice $C_n^u$ satisfying (\ref{Cn}) (which is fixed under the bootstrap distribution).
 \end{theorem}
 The proof of Theorem \ref{thnpbootmle} is presented in  Appendix \ref{AppendixB}. 
 In Appendix B we first establish that $d_{n1}(Q_n^{\#},Q_n)=O_P(n^{-1/2-\alpha/4})$, and we use that, in combination with $d_{10}(Q_n,Q_0)=O_P(n^{-1/2-\alpha/4})$, this also implies $d_{01}(Q_n^{\#},Q_0)=O_P(n^{-1/2-\alpha/4})$. Thus, clearly, $d_{n1}(Q_n^{\#},Q_n)$ is an equally powerful dissimilarity as $d_{01}()$.
In fact, assuming that ${\cal M}$ has the extra model structure (\ref{calFmodel}), in Appendix D  we also explicitly show that $d_{n1}(Q_n^{\#},Q_n)$ dominates a specified quadratic dissimilarity. 


 Note that if $C^u=C_n^u$, then conditional on $P_n$, $C_n^u$ is still fixed, so that establishing the last result  in Theorem \ref{thnpbootmle} only requires checking that the proof of the stated convergence of the bootstrapped HAL-MLE $(Q_{n,C_Q^u}^{\#},G_{n,C_G^u}^{\#})$ to the HAL-MLE $(Q_{n,C_Q^u},G_{n,C_G^u})$ at  a fixed $C^u=(C_Q^u,C_G^u)$  w.r.t. the loss-based dissimilarities $d_{n1}$ and $d_{n2}$ holds uniformly in $C^u$ between 
 the true sectional variation norms $C_0^u$ and the model upper bound $C^u$. The validity of this result does not even rely on  $C_n^u$  exceeding $C_0^u$, but the latter is needed for establishing that the HAL-MLE $Q_{n,C_n^u}$ is consistent for $Q_0$ and thus the efficiency of the HAL-TMLE $\Psi(Q_n^*)$.



\subsection{Preservation of rate of convergence for the targeted bootstrap estimator}
In Appendix \ref{AppendixC} we prove that $d_{01}(Q_n^{\#*},Q_0)=O_P(n^{-1/2-\alpha/4})$, under the same conditions as assumed in our general Theorem \ref{thnpboothaltmle}.


\subsection{The nonparametric bootstrap for the  HAL-TMLE}
We can now imitate the efficiency proof for the HAL-TMLE to obtain the desired result for the bootstrapped HAL-TMLE of $\Psi(Q_n^*)$.
By Theorem \ref{thnpbootmle}, under the assumptions of Theorem \ref{thefftmle} for asymptotic efficiency of the TMLE, 
we have that 
all five terms $d_{n1}(Q_n^{\#},Q_n)$, $d_{01}(Q_n^{\#},Q_0)$, $d_{01}(Q_n^{\#*},Q_0)$, $d_{n2}(G_n^{\#},G_n)$, $d_{02}(G_n^{\#},G_0)$ are $O_P(n^{-1/2-\alpha/4})$.
For a model with extra structure (\ref{calFmodel}), we consider the bootstrap for a data adaptive selector $C_n^u=(C_{Qn}^u,C_{Gn}^u)$ satisfying
(\ref{Cn}). A general model ${\cal M}$ might also be indexed by a universal bound $C$ for some quantity $C(P)$ for any $P\in {\cal M}$, which could then also be data adaptively selected as long as it satisfies (\ref{Cn}) with $C_0=C(P_0)$.

\begin{theorem}\label{thnpboothaltmle}

{\bf Assumptions:}
Consider the statistical model ${\cal M}$ and target parameter $\Psi:{\cal M}\rightarrow\openr$ s
satisfying
(\ref{M20bounds}), (\ref{sectionalvarbound}), (\ref{boundingR2}), (\ref{contDstar}).
Consider the above defined HAL-MLEs $Q_n$, $G_n$ satisfying, with probability tending to 1,
$P_n L_1(Q_n)\leq P_n L_1(Q_0)$ and $P_n L_2(G_n)\leq P_n L_2(G_0)$. Consider also the above defined bootstrapped HAL-MLEs $Q_n^{\#}$, $G_n^{\#}$ satisfying, with probability tending to 1, conditional on $(P_n:n\geq 1)$,  $P_n^{\#}L_1(Q_n^{\#})\leq P_n L_1(Q_n)$ and $P_n^{\#}L_2(G_n^{\#})\leq P_n^{\#}L_2(G_n)$. 
Consider the HAL-TMLE $Q_n^{\#*}=Q_{n,\epsilon_n^{\#}}^{\#}$ and assume (\ref{efficeqnb}) $r_n^{\#}=P_n^{\#}D^*(Q_n^{\#*},G_n^{\#})=o_P(n^{-1/2})$.

{\bf TMLE is efficient:}
The standardized TMLE is asymptotically efficient:  $Z_n^1\equiv n^{1/2}(\Psi(Q_n^*)-\Psi(Q_0))\Rightarrow_ d N(0,\sigma^2_0)$, where $\sigma^2_0=P_0D^*(Q_0,G_0)^2$.

{\bf Bootstrapped HAL-MLE:}
$d_{01}(Q_n^{\#},Q_n)=O_P(n^{-1/2-\alpha/4})$, $d_{02}(G_n^{\#},G_0)=O_P(n^{-1/2-\alpha/4})$ and $d_{01}(Q_n^{\#*},Q_0)=O_P(n^{-1/2-\alpha/4})$.

{\bf Bootstrapped HAL-TMLE:}
 Conditional on $(P_n:n\geq 1)$, the bootstrapped TMLE is asymptotically linear:
\[
\Psi(Q_n^{\#*})-\Psi(Q_n)=(P_n^{\#}-P_n)D^*(Q_n,G_n)+O_P(n^{-1/2-\alpha/4}).\]
As a consequence, conditional on $(P_n:n\geq 1)$, the standardized bootstrapped TMLE converges to $N(0,\sigma^2_0)$: $Z_n^{1,\#}\equiv n^{1/2}(\Psi(Q_n^{\#*})-\Psi(Q_n^*))\Rightarrow_d N(0,\sigma^2_0)$.
\newline
{\bf Consistency of the nonparametric bootstrap for HAL-TMLE at data adaptive selector $C_n^u$:}
Assume the extra model structure (\ref{calFmodel}) on ${\cal M}$, and its corresponding definitions of the HAL-MLEs indexed by sectional variation norm bounds $C=(C^u,C^l)$.
This theorem can be applied to the bootstrap distribution at a data adaptive $C_n=(C_n^u,C_n^l)$ satisfying (\ref{Cn}).
\end{theorem}
The proof of this theorem is presented in Appendix \ref{AppendixD1}.

\section{Finite sample modifications of the nonparametric bootstrap distribution for model with extra structure (\ref{calFmodel})}
\label{sectfour} \

In this section we focus on the case that the model ${\cal M}$ satisfies the extra structure (\ref{calFmodel}). The finite sample modifications proposed here are evaluated in our simulation study in Section \ref{sectsix} for our two examples. 
The nuisance parameter estimates $Q_n$ and $G_n$ are key inputs of the HAL-TMLE
bootstrap. The HAL estimations of these nuisance parameters depend largely on
the selection of the upper bound of the sectional variation norm $C^u=(C_Q^u,C_G^u)$. 
We will focus on a data adaptive selector of $C_{Qn}^u$ (replacing $C_Q^u$), for a given selector $C_{Gn}^u$, where the latter is chosen to be the cross-validation selector. Since our target parameter is a function of $Q$ only, we suggest that the selection of $C_{Qn}^u$ is fundamentally more important than $C_{Gn}^u$, and also creates enough room for our desired finite sample adjustment of the nonparametric bootstrap.
In the software implementation of LASSO, the $L_1$-norm constraint $C_Q^u$ is translated into a penalized empirical risk with $L_1$-penalty hyper-parameter $\lambda$, where a choice of $C_Q^u$ corresponds with a unique choice $\lambda$. 
In the sequel, we will propose a selector of $\lambda$, and thereby of $C_Q^u$.

Ideally, we want to set $C_Q^u=C_{Q0}^u$ equal to the sectional variation norm of $Q_0$, so that the bootstrap model for the HAL-MLE $Q_n^{\#}$  is large enough for unbiased estimation of $Q_n$. Due to the asymptotic equivalence of the cross-validation selector $C_{Qn,CV}^u$ with the oracle selector that optimizes the loss-based dissimilarity, the cross-validation selector $C_{Qn,CV}^u$ will approximate $C_{Q0}^u$ as sample size increases. However, in finite samples, when the true sectional variation norm $C_{Q0}^u$
of $Q_0$ is large ($\lambda_0$ is small), the cross-validation selector $C_{Qn,CV}^u$ will 
 tend to be smaller than the oracle
value $C_{Q0}^u$ ($\lambda_{CV} > \lambda_0$),  That is, $C_{Qn,CV}^u$ optimally trades off bias and variance for estimation of $Q_0$, but fixing $C_Q^u$ at this choice $C_{Qn,CV}^u$ might oversimplify the complexity of the target $Q_n^*$ of the bootstrap distribution, and thereby causes the bootstrap to  under-estimate the variability of the true sampling distribution of the TMLE.
As a result, the bootstrap confidence interval will potentially still be anti-conservative.

Since the oracle choice $\lambda_0$ is unknown, we propose to estimate $\lambda_0$ with a plateau selection method. 
Consider a pre-specified ordered (from large to small) sequence of lambda candidates $\Lambda = (\lambda_1, \lambda_2, ..., \lambda_J)$ with corresponding HAL-MLEs $Q_{n,\lambda_j}$ and HAL-TMLEs $Q_{n,\lambda_j}^*$, $j=1,\ldots,J$. We set $\lambda_1=\lambda_{n,CV}$ so that we only consider sectional variation norm constraints larger than the cross-validation selector $C_{Qn,CV}^u$. The sectional variation norm of $Q_{n,\lambda_j}$ will thus be increasing in $j$. For each $\lambda_j$ we compute the width $w_j=(q_{0.975,n,\lambda_j}^{\#}-q_{0.025,n,\lambda_j}^{\#})\sigma_n$ of the nonparametric bootstrap confidence interval based on bootstrapping the standarized TMLE $n^{1/2}(\Psi(Q_{n,\lambda_j}^*)-\Psi(Q_0) ) / \sigma_n$, 
given by $[\Psi(Q_{n}^*)+q_{0.025,n,\lambda_j}^{\#}\sigma_n,\Psi(Q_n^*)+q_{0.975,n,\lambda_j}^{\#}\sigma_n]$, $j=1,\ldots,J$.
The interval widths monotonically increase and should generally show de-acceleration around $\lambda_0$ where it will move towards a plateau, and, eventually it might become erratic.
A related theoretical reference of this phenomena is Theorem 1 in \cite{Davies&vanderLaan14}, which proposes such a plateau selector, and proves the consistency of the corresponding plateau variance estimator in a growing model that eventually captures the true data distribution (even though, in practice, plateau appear in much greater generality). As the model grows, the variance estimator keeps increasing, but once the model contains the true distribution the variance estimator is consistent/unbiased. Therefore, for large sample sizes we should see the same or similar variance estimate, but once model gets too big relative to sample size (defined in the Theorem 1 in \cite{Davies&vanderLaan14}), it becomes too variable and erratic. In our methodology, the widths of the confidence intervals correspond with variance/uncertainty estimation, and our models grow due to increase in sectional variation norm, and, indeed, as the sectional variation norm passes the variation norm of the true function, the growing model captures the truth.

 So a similar intuition holds for our estimator.
If we set variation norm $C^u_Q$ smaller than true $C^u_{Q_0}$, and let $n$ go to infinity, we are inconsistent (negatively biased) for the true variance. If we set $C^u_Q > C^u_{Q_0}$, any TMLE is efficient so will have the same asymptotic variance. In between as we increase $C^u_Q$ towards $C^u_{Q_0}$, the width of the confidence intervals grows accordingly. Therefore, we expect for large sample size to see that the width curve will increase as $C^u_Q$ moves towards $C^u_{Q_0}$ and become flat after $C^u_Q > C^u_{Q_0}$. 
Through numerical simulations, we indeed observed that $\lambda_0$ is near where the plateau begins. It remains to decide on a method for determining the location of the start of the de-acceleration. A variety of methods could be proposed here. In our concrete implementation demonstrated in our simulation study, we compute the location of the start of the plateau as the location at which the second derivative is maximized, where we use the $\log \lambda$-scale (due to $\lambda$ having very small values). Specifically, 
$\lambda_{plateau} \triangleq \lambda_j$, where 
\begin{align}
j = \argmax_{j = 2,...,J - 1}
\frac{(w_{j + 1} - w_j) - (w_j - w_{j - 1})}{(\log ({\lambda _{j + 1}}) - \log
(\lambda _j))(\log (\lambda _j) - \log (\lambda _{j - 1}))}    
\label{eq:2_derivative}
\end{align}
We choose a log-uniform grid of pre-specified $\lambda$ to simplify the finite difference estimation of the derivative, and we leave it an important future work to implement a potentially better estimator with more flexible choice of $\lambda$ grid.

\begin{figure}[!ht]
    \centering
    \includegraphics[width=0.95\textwidth]{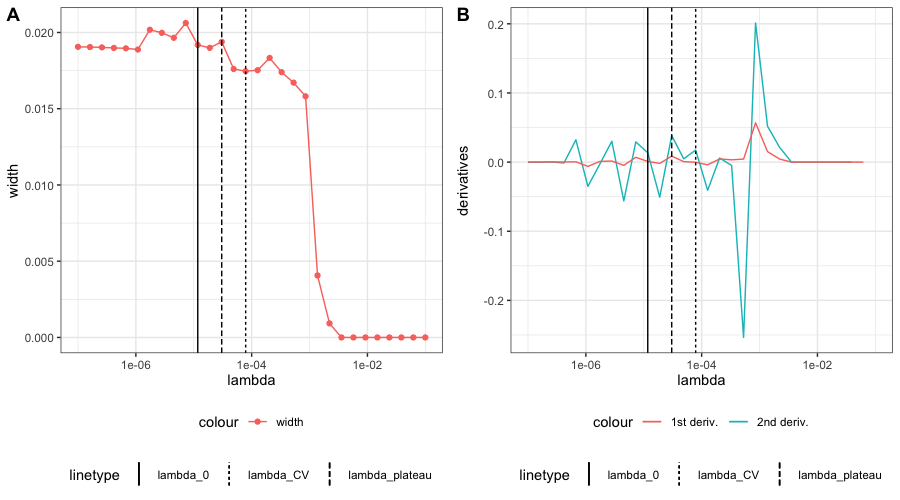}
    \caption{(A) A simulated example of Wald-type interval width as a function of $\lambda$. (B) The first and second order derivatives of the same curve. Vertical lines indicate $\lambda_0$, $\lambda_{CV}$ and $\lambda_{plateau}$.}
    \label{fig:plateau_exist}
\end{figure}

Figure \ref{fig:plateau_exist} illustrates a simulated example of the curve $\log(\lambda)\rightarrow w(\lambda)$. As the value of $\lambda$ decreases starting at $\lambda_{CV}$, we observe a slow increase initially (almost a flat area around $\lambda_{CV}$), then an accelerated increase, till it starts reaching its plateau right after $\lambda_0$. Our method looks for the numerical maximum of (discrete) second-order derivative (\ref{eq:2_derivative}), where the function starts moving towards the plateau. Another method might be to look for the actual start of the plateau, but our concern is that this might corresponds with a plateau due to pure overfitting the data (where the finite sample only allows so much overfitting).

\paragraph{Increasing the scaling $\sigma_n$-factor by taking into account bias of bootstrap sampling distribution}\label{sec:rmse_scaled}
Another modification we propose concerns the bias of the bootstrap distribution. We assume that we used the above method for selecting a $\lambda_n = \lambda_{plateau}$. We will use as point estimate $\Psi(Q_n^*)$, where $Q_n^* = Q_{n,\lambda_{n,CV}}^*$, i.e, the TMLE using the cross-validated HAL-MLE. So the role of the bootstrap is to determine a confidence interval around this point estimate. Our confidence interval will be of the form $[\Psi(Q_n^*) + q_{n,0.025}^{\#}\sigma_n^{\#}/n^{1/2}, \Psi(Q_n^*) + q_{n,0.975}^{\#}\sigma_n^{\#}/n^{1/2}]$, where we use the nonparametric bootstrap at fixed sectional variation norm implied by $\lambda_n$, but centered to have mean zero, to obtain these two quantiles. The bias in the bootstrap distribution will instead be incorporated in  $\sigma_n^{\#}$  by defining $\sigma_n^{\#2}$ as the MSE of the bootstrap realizations $\Psi(Q_{n,i}^{\#*})$ relative to $\Psi(Q_n^*)$, $i=1, \ldots, N$, where $N$ is the number of bootstrap samples drawn from $P_n$.

The motivation is that in general the nonparametric bootstrap will also inherit bias of the sampling distribution of $n^{1/2}(\Psi(Q_n^*)-\Psi(Q_0))/\sigma_n$. For example, if there is finite sample bias of $\Psi(Q_n^*)$ that is hurting the coverage of a Wald-type confidence interval, the bootstrap distribution (i.e., its quantiles) will likely further bias in the same direction. We choose not to estimate the bias with the bootstrap and compensate the bootstrap distribution accordingly through shifting it, since estimates of bias are typically unreliable.
Instead, we  widen the bootstrap confidence interval by replacing the scaling factor $\sigma_n$ by the square root of the MSE of $\Psi(Q_n^{\#*})$ w.r.t. $\Psi(Q_n^*)$.
Specifically, the ``RMSE-scaled bootstrap'' takes the form 
\begin{align}
  [\Psi(Q_n^*)+ \sigma^\#_n q^{\#}_{n,0.025}/n^{1/2}, \Psi(Q_n^*) + \sigma^\#_n q^{\#}_{n,0.975}/n^{1/2}],
\end{align}
where (using short-hand notation) 
\[
\sigma^\#_n \triangleq \sqrt{ \frac{1}{N} \sum_{i =
1}^N{(\Psi^{\#*}_{i,n} - \Psi(Q_n^*)) ^ 2} } = \sqrt{\mbox{bias}(\Psi^{\#*}_{i,n})^2 +
\mbox{stddev}(\Psi^{\#*}_{i,n})^2}\] is the estimated RMSE of the bootstrap estimator
$\Psi^{\#*}_{i,n}=\Psi(Q_{n,i}^{\#*})$, and $q^{\#}_{n,\alpha}$ is the $\alpha$-quantile of the
bootstrap distribution of standardized  $Z^{\#}_{i,n} =
n^{1/2}(\Psi^{\#*}_{i,n} - \frac{1}{N} \sum_{i = 1}^N \Psi^{\#*}_{i,n}) / \mbox{stddev}(\Psi^{\#*}_{i,n})$. 

The full modified HAL-TMLE bootstrap procedure we propose in this article can be
summarized in the following pseudo-algorithm:

\begin{algorithm}[H]
\caption{modified HAL-TMLE bootstrap procedure}
pre-specify a grid of $\lambda$ values, $\Lambda$\;
\For{$\lambda \in \Lambda$}{
    fit HAL-MLE $Q_n$ using tuning parameter $\lambda$\;
    perform HAL-TMLE and record point TMLE $\Psi^*_n(\lambda)$\;
}
perform cross-validation to select $\lambda_{CV}$; record the HAL-TMLE point estimate $\Psi(Q^*_n)$ with $Q_n^*=Q_{n,\lambda_{CV}}^*$\;
Compute the plateau selector  $\lambda_{plateau}$ among $\lambda \geq \lambda_{CV}$  based on running the nonparametric bootstrap for $n^{1/2}(\Psi(Q_{n,\lambda}^{\#*})-\Psi(Q_{n,\lambda}^*) )/\sigma_{n,\lambda}^{\#}$\;
Set $\lambda=\lambda_{plateau}$, perform HAL-TMLE bootstrap $N$ times to obtain quantiles $q_{n,0.025}^{\#},q_{n,0.975}^{\#}$ of $n^{1/2}(\Psi(Q_{n,\lambda}^{\#*})-E_{P_n}\Psi(Q_{n,\lambda}^{\#*}))/\sigma_{n,\lambda}^{\#}$\;
 compute $\sigma^\#_n = \sqrt{ \frac{1}{N} \sum_{i = 1}^N{(\Psi(Q^{\#^*}_{i,n}) - \Psi(Q^*_n)) ^ 2}}$ \;
report $\Psi(Q^*_n)$ as the final point estimator; report the 95\% confidence interval of the target parameter as $[\Psi(Q^*_n) + \sigma^\#_n q^{\#}_{n,0.025}/n^{1/2}, \Psi(Q^*_n) + \sigma^\#_n q^{\#}_{n,0.975}/n^{1/2}]$.
\end{algorithm}

\section{Examples}\label{sectfive}
In this section we apply our general theorem, by verifying its conditions, for asymptotic consistency of the nonparametric bootstrap of HAL-TMLE to two examples involving a nonparametric model. 
In the next section we will actually implement our nonparametric bootstrap based confidence intervals for these two examples, carry out a simulation study, and evaluate its practical performance w.r.t. finite sample coverage. 

\subsection{Nonparametric estimation of average treatment effect}
Let $O=(W,A,Y)\sim P_0$, where $W\in [0,\tau_1]\subset \openr^{m_1}_{\geq 0}$ is an $m_1$-dimensional vector of baseline covariates, $A\in \{0,1\}$ is a binary treatment, and $Y\in \{0,1\}$ is a binary outcome. For a possible data distribution $P$, let $\bar{Q}(P)=E_P(Y\mid A,W)$, $
\bar{G}(P)=P(A=1\mid W)$, and let $Q_W(P)$ be the cumulative  probability distribution of $W$. 
Let $Q_1=Q_W$, $Q_2=\mbox{logit}\bar{Q}$, $Q=(Q_1,Q_2)$, and $G=\mbox{logit}\bar{G}$.  Let $g(a\mid W)=P(A=a\mid W)=\bar{G}(W)^a(1-\bar{G}(W))^{1-a}$.
In addition, let $m_{11}=m_1$ and $m_{12}=m_1+1$, in terms of our general notation. Suppose that our model assumes that $\bar{G}(W)$ depends on a  possible subvector of $W$, and let $m_2 $ be the dimension of this subvector.

\paragraph{Statistical model:}
Since $Q_1=Q_W$ is a cumulative distribution function, it is a monotone $m_1$-variate cadlag function and its sectional variation norm equals its total variation which thus equals 1. 
We assume that $Q_2$ is an element of the class of $m_{12}$-dimensional cadlag functions with sectional variation norm bounded from above by some $C_{Q2}^u$. Here one can treat $A$ as continuous on $[0,1]$ and assume that $Q_2$ is a step-function in $A$ with single jump at 1, allowing us to embed functions of continuous and discrete covariates in a cadlag function space.
Similarly, we assume $G$ is an element of the class of $m_2$-dimensional cadlag functions with sectional variation norm bounded by a $C_G^u$.
Let's denote these parameter spaces for $Q_1,Q_2$ and $G$ with ${\cal F}_{11}$, ${\cal F}_{12}$ and ${\cal F}_2$, respectively.
Let ${\cal F}_1={\cal F}_{11}\times {\cal F}_{12}$ be the parameter space of $Q=(Q_1,Q_2)$.
For a given $C_Q^u=(C_{Q1}^u=1,C_{Q2}^u),C_G^u<\infty$, consider the statistical model \begin{equation}\label{modelate}
{\cal M}=\{P: Q_W\in {\cal F}_{11}, \bar{Q}\in {\cal F}_{12}, G\in {\cal F}_2\}.\end{equation}
 Thus, ${\cal M}$ is defined as the set of all possible probability distributions for which the logit of the conditional means of $Y$ and $A$ are cadlag functions with sectional variation norm bounded by $C_Q^u$ and $C_G^u$, respectively. 
 Since $\mbox{logit}\bar{G}$ and $\mbox{logit}\bar{Q}$ are bounded in supremum norm (implied by their bounds on the sectional variation norm), it follows that $\bar{G}$ and $\bar{Q}$ are bounded from below by $\delta>0$ and from above by $1-\delta$ for some $\delta>0$. We will refer to this bound $\delta=\delta(C_Q^u,C_G^u)$ separately in our bounds below, even though it is implied by the sectional variation norm bound $C^u$.
 In particular, this implies the strong positivity assumption $\min_{a\in \{0,1\}}g(a\mid W)>\delta>0$ $Q_W$-a.e. 
 
 Notice that indeed our parameter space for $Q=(Q_1,Q_2)$ and of $G$ is of type (\ref{calFmodel}) or (\ref{calFmodelplus}). Specifically, $Q_W$ is of the type (\ref{calFmodelplus}) with $C_{Q1}^l=C_{Q1}^u=1$, while $Q_2$ and $G$ have parameter spaces of type (\ref{calFmodel}) with only an upper bound $(C_{Q2}^u,C_G^u)$ on their sectional variation norm.  This demonstrates that our model ${\cal M}$ is represented as the general  model formulation defined in Section 2.
 
\paragraph{Target parameter:}
Let $\Psi:{\cal M}\rightarrow\openr$ be defined by $\Psi(P)=\Psi_1(P)-\Psi_0(P)$, where $\Psi_a(P)=E_PE_P(Y\mid A=a,W)$.
Note that $\Psi(P)$ only depends on $P$ through $Q(P)=(Q_1,Q_2)$, so that we will also use the notation $\Psi(Q)$ instead of $\Psi(P)$.
Let's focus on $\Psi_1(P)$ which will also imply the formulas for $\Psi_0(P)$ and thereby $\Psi(P)$.

\paragraph{Loss functions for $Q$ and $G$:}
Let $L_{11}(Q_W)=\int_x(I(W\leq x)-Q_W(x))^2r(x) dx$ for some weight function $r>0$  be the loss function for $Q_{10}=Q_{W,0}$. 
Let 
\begin{align*}
d_{01}(Q_W,Q_{W,0})=P_0L_{11}(Q_W)-P_0L_{11}(Q_{W,0})    
\end{align*}
be the corresponding  loss-based dissimilarity.
Let $L_{12}(Q_2)=-\{Y\log\bar{Q}(A,W)+(1-Y)\log(1-\bar{Q}(A,W))\}$ be the log-likelihood loss function for the conditional mean $\bar{Q}_0$ and thereby $Q_{20}=\mbox{logit}\bar{Q}_0$. Let
$d_{02}(Q_2,Q_{20})=P_0 L_{12}(Q_2)-P_0L_{12}(Q_{20})$ be the corresponding Kullback-Leibler dissimilarity.
We can then define the sum-loss $L_1(Q)=L_{11}(Q_1)+L_{12}(Q_2)$ for $Q_0=(Q_{10},Q_{20})$, and its loss-based dissimilarity
 $d_{01}(Q,Q_0)=P_0L_1(Q)-P_0L_1(Q_0)$ which equals the sum of the following two dissimilarities:
 \begin{eqnarray*}
d_{01}(Q_1,Q_{10})&=&\int_x (Q_{W}(x)-Q_{W,0}(x))^2 r(x) dx, \\
d_{02}(Q_2,Q_{20})&=&\int \log \left[\left(\frac{\bar{Q}_0}{\bar{Q}}\right)^y \left( \frac{1-\bar{Q}_0}{1-\bar{Q}}\right)^{1-y}\right](a,w) 
 dP_0(w,a,y) ,\end{eqnarray*}
 $\bar{Q}=\bar{Q}(Q_2)$ and $\bar{Q}_0=\bar{Q}(Q_{20})$ are implied by $Q_2$ and $Q_{20}$, respectively.
Let $L_2(G)=-\{A\log \bar{G}(W)+(1-A)\log(1-\bar{G}(W))\}$ be the loss function for $G_0=\mbox{logit}P_0(A=1\mid W)$, and let $d_{02}(G,G_0)=P_0L_2(G)-P_0L_2(G_0)$ be the Kullback-Leibler dissimilarity between $G$ and $G_0$.

\paragraph{Canonical gradient and corresponding exact second order expansion:}
The canonical gradient of $\Psi_a$ at $P$ is given by:
\[
D^*_a(Q,G)=\frac{I(A=a)}{g(A\mid W)}(Y-\bar{Q}(A,W))+\bar{Q}(a,W)-\Psi_a(Q).\]
The exact second-order remainder $R_{20}^a(P,P_0)\equiv \Psi_a(P)-\Psi_a(P_0)+P_0 D^*_a(P)$ is given by:
\[
R_{20}^a(Q_2,G,Q_{20},G_0) =\int \frac{(g-g_0)(a\mid w)}{g(a\mid w)}(\bar{Q}-\bar{Q}_0)(a,w) dP_0(w).\]

\paragraph{Bounding the second order remainder:}
By using Cauchy-Schwarz inequality, we obtain the following bound on $R_{20}^a(P,P_0)$:
\[
| R_{20}^a(P,P_0)| \leq \delta^{-1}\pl \bar{Q}_a-\bar{Q}_{a0}\pl_{P_0}\pl G-G_0\pl_{P_0} ,\]
where $\bar{Q}_a(W)=\bar{Q}(a,W)$, $a\in \{0,1\}$.
Thus, $D^*(P)=D^*_1(P)-D^*_0(P)$, $R_{20}(P,P_0)=R_{20}^1(P,P_0)-R_{20}^0(P,P_0)$, and the upper bound for $R_{20}(P,P_0)$ can be defined as the sum of the two upper bounds for $R_{20}^a(P,P_0)$ in the above inequality, $a\in \{0,1\}$.

By \cite{vanderVaart98} we have $\pl p^{1/2}-p_0^{1/2}\pl_{P_0}^2\leq  P_0 \log p_0/p $, where $p$ and $p_0$ are densities of $P$ and $P_0$, with $P_0 \ll P$. For Bernoulli distributions,  we have
$\pl p-p_0\pl^2_{P_0}\leq 4 \pl p^{1/2}-p_0^{1/2}\pl^2_{P_0}\leq 4 P_0\log p_0/p$. 
Following the same proof as in Lemma 4 of \cite{vanderLaan15}, we note $p$ is playing role of $p(Y|A,W)$. so $d_0(p,p_0)$ is our KL dissimilarity $d_{02}(\bar{Q},\bar{Q}_0)$. It now remains to show $\|p-p_0\|^2 = \int_{a,w} \int_y (p-p_0)^2(y|a,w) dP_0(y|a,w) dP_0(a,w)$, where $\int_y$ is just sum over $y=0$ and $y=1$. Since $Y$ is binary, $p(y=0|a,w) = 1 - p(y=1|a,w) = 1 - \bar{Q}(a, w)$, and we obtain  $\int (\bar{Q}-\bar{Q}_0)^2(a,w)dP_0(a,w)\leq 4 d_{02}(\bar{Q},\bar{Q}_0)$ and thus 
$\pl \bar{Q}_a-\bar{Q}_{a0}\pl^2_{P_0}\leq 4\delta^{-1}d_{02}(\bar{Q},\bar{Q}_0)$.
Therefore, we conclude that $\pl \bar{Q}_a-\bar{Q}_{a0}\pl_{P_0}\leq 2\delta^{-1/2} d_{02}^{1/2}(\bar{Q},\bar{Q}_0)$.
Similarly, it follows that $\pl G-G_0\pl_{P_0}\leq 2 d_{02}^{1/2}(G,G_0)$.
This thus shows the following bound on $R_{20}^a(P,P_0)$:
\[
| R_{20}^a(P,P_0)| \leq 2\delta^{-1.5} d_{02}^{1/2}(\bar{Q},\bar{Q}_0) d_{02}^{1/2}(G,G_0).\]
The right-hand side represents the function $f({\bf d}_{01}^{1/2}(Q,Q_0),{\bf d}_{02}^{1/2}(G,G_0))$ for the parameter $\Psi_a$  in our general notation: $f(x=(x_1,x_2),y)=
4 \delta^{-1.5} x_2 y$.
The sum of these two bounds for $a\in \{0,1\}$ (i.e, $2f()$) provides now a conservative bound for  $R_{20}=R_{20}^1-R_{20}^0$:
\begin{equation}\label{r2upperboundexample1}
| R_{20}(P,P_0)| \leq f(d_{02}^{1/2}(\bar{Q},\bar{Q}_0),d_{02}^{1/2}(G,G_0))\equiv 8\delta^{-1.5} d_{02}^{1/2}(\bar{Q},\bar{Q}_0) d_{02}^{1/2}(G,G_0).\end{equation}
This verifies (\ref{boundingR2}).
We note that this bound is very conservative due to the arguments we provided in general in the previous section for double robust estimation problems.

\paragraph{Continuity of canonical gradient:}
Regarding the continuity assumption (\ref{contDstar}), we note that $P_0\{D^*_a(P)-D^*_a(P_0))^2$ can be bounded by $\pl G-G_0\pl_{P_0}^2+\pl \bar{Q}_a-\bar{Q}_{a0}\pl^2_{P_0}$ and $(\Psi_a(Q)-\Psi_a(Q_0))^2$, where the constant depends on $\delta$. The latter square difference can be bounded in terms of $\pl \bar{Q}_a-\bar{Q}_{a0}\pl^2_{P_0}$ and  by applying our integration by parts formula to  $\int \bar{Q}_a(w) d(Q_W-Q_{W0})(w)$ by $d_{01}(Q_W,Q_{W0})$, where the multiplicative constant depends on $C_Q^u$. 
We conclude that $P_0\{D^*_a(P)-D^*_a(P_0))^2$ is bounded in terms of $d_{01}(Q,Q_0)+d_{02}(G,G_0)$.
Thus this proves (\ref{contDstar}) for $D^*=D^*_1-D^*_0$.

\paragraph{Uniform model bounds on sectional variation norm:} It also follows immediately that the sectional variation norm model bounds $M_1,M_2,M_3$ (\ref{sectionalvarbound})
of $L_1(Q)$, $L_2(G)$ and $D^*(P)$ are all finite, and can be expressed in terms of $(C_Q^u,C_G^u,\delta)$. This verifies the model assumptions of Section 2.


\paragraph{HAL-MLEs:}
 Let $Q_n=\argmin_{Q\in {\cal F}_1}P_n {\bf L}_1(Q)$ and $G_n=\argmin_{G\in {\cal F}_2}P_n L_2(G)$ be the HAL-MLEs. 
 As shown in \citep{vanderLaan15,Benkeser&vanderLaan16}, 
 $\bar{Q}_n$ and $G_n$ can be computed with standard LASSO logisitic regression software using a linear logistic regression model with around $n 2^{m_1}$ indicator basis functions, where $m_1$ is the dimension of $W$.

 Note that $Q_{W,n}$ is just an unrestricted MLE and thus equals the empirical cumulative distribution function. 
 Therefore, we actually have that $\pl Q_{W,n}-Q_{W,0}\pl_{\infty}=O_P(n^{-1/2})$ in supremum norm, while $d_{02}(Q_{2n},Q_{20})$ and $d_{02}(G_n,G_0)=O_P(n^{-1/2-\alpha/4})$
 where $d$ is the dimension of $O$. If $m_2<d-2$, then one  should be able to improve the bound into $n^{-1/2-\alpha(m_2)}$.
 
\paragraph{CV-HAL-MLEs:}
 The above HAL-MLEs are determined by $(C_Q^u=(1,C_{Q2}^u),C_G^u)$ and could thus be denoted with $Q_{n,C_Q^u}=\hat{Q}_{C_Q^u}(P_n)$ and $G_{n,C_G^u}=\hat{G}_{C_G^u}(P_n)$.
 Let $C_{Q0}=\pl Q_0\pl_v^*=(1,\pl Q_{20}\pl_v^*)$ and $C_{G0}=\pl G_0\pl_v^*$, respectively, which are thus smaller than $C_Q^u$ and $C_G^u$, respectively.
 We can now define the cross-validation selector that selects the best HAL-MLE over all $C_{Q}$ and $C_G $ smaller than these upper-bounds:
 \begin{eqnarray*}
 C_{Qn}&=&\argmin_{C_{Q1}=1,C_{Q2}<C_{Q2}^u}E_{B_n}P_{n,B_n}^1L_1(\hat{Q}_{C_Q}(P_{n,B_n}^0)) \\
C_{Gn}&=&\argmin_{C_G<C_G^u}E_{B_n}P_{n,B_n}^1L_2(\hat{G}_{C_G}(P_{n,B_n}^0)),
\end{eqnarray*}
where $B_n\in \{0,1\}^n$ is a random split in training sample $\{O_i:B_n(i)=0\}$ with empirical measure $P_{n,B_n}^0$ and validation sample
$\{O_i:B_n(i)=1\}$ with empirical measure $P_{n,B_n}^1$.
This defines now the CV-HAL-MLE $Q_n=Q_{n,C_{Qn}}$ and $G_n=G_{n,C_{Gn}}$ as well. Thus, by setting $C_Q^u=C_{Qn}$ and $C_G^u=C_{Gn}$, our HAL-MLEs equal the CV-HAL-MLE.

\paragraph{HAL-TMLE:}
 Let $\mbox{logit}\bar{Q}_{n,\epsilon}=\mbox{logit}\bar{Q}_n+\epsilon C(G_n)$, or, equivalently, $Q_{2n,\epsilon}=Q_{2n}+\epsilon C(G_n)$, where $C(G_n)(A,W)=(2A-1)/g_n(A\mid W)$.
 Let $\epsilon_n=\argmin_{\epsilon}P_n L_{11}(\bar{Q}_{n,\epsilon})$. This defines the TMLE $\bar{Q}_n^*=\bar{Q}_{n,\epsilon_n}$ of $\bar{Q}_0$, and thereby $Q_{2n}^*=Q_{2n,\epsilon_n}$. We can also define a local least favorable submodel $\{Q_{W,n,\epsilon_2}:\epsilon_2\}$ for $Q_{W,n}$ but since $Q_{W,n}$ is an NPMLE one will have that $\epsilon_{2n}=\argmin_{\epsilon_2}P_n L_{11}(Q_{W,n,\epsilon_2})=0$, and thereby that the TMLE of $Q_0$ for any such 2-dimensional least favorable submodel is given by  
 $Q_n^*=(Q_{W,n},{Q}_{2n}^*)$. It follows that $P_n D^*(Q_n^*,G_n)=0$.
 
\paragraph{Preservation of rate for HAL-TMLE:}
 Lemma \ref{lemmathalmle} in Appendix \ref{AppendixA}
shows $d_{01}(Q_n^*,Q_0)$  converges at same rate $O_P(n^{-1/2-\alpha/4})$ as $d_{01}(Q_n,Q_0)$.

\paragraph{Asymptotic efficiency of HAL-TMLE and CV-HAL-TMLE:}
 Application of Theorem \ref{thefftmle} shows that $\Psi(Q_n^*)$ is asymptotically efficient, where one can either choose $Q_n$ as  a fixed HAL-MLE using $C_Q=C_Q^u$ or the CV-HAL-MLE using $C_Q=C_{Qn}$, and similarly, for $G_n$. The preferred estimator would be the CV-HAL-TMLE.
 
\paragraph{Asymptotic validity of the nonparametric bootstrap for the HAL-MLEs:}
 Firstly, note that the bootstrapped HAL-MLEs \[
 {Q}_{2n}^{\#}=\argmin_{\pl Q_2\pl_v^*< C_{Q2}^u,Q_2\ll {Q}_{2n}}P_n^{\#}L_{12}(Q_2),\]
  and 
 $$G_n^{\#}=\argmin_{\pl G\pl_v^*<C_G^u,G\ll G_n}P_n^{\#}L_2(G)$$ are easily computed as a standard LASSO regression using $L_1$-penalty $C_{Q2}^u$ and $C_G^u$ and including the $\approx n$ indicator basis functions with the non-zero coefficients selected by $Q_n$ and $G_n$, respectively. This makes the actual computation of the nonparametric bootstrap distribution a very doable computational problem, even though the single computation of $Q_n$  and $G_n$ is highly demanding for large dimension of $W$ and sample size $n$.
 $Q_{1n}^{\#}=Q_{W,n}^{\#}$ is simply the empirical probability measure of $W_{1}^{\#},\ldots,W_n^{\#}$ by sampling $n$ i.i.d. observations from $Q_{W,n}$.

\paragraph{Behavior of HAL-MLE under sampling from $P_n$:}
By Theorem \ref{thnpbootmle} we have that each of the terms $d_{n12}(Q_{2n}^{\#},Q_{2n})$;
 $d_{n2}(G_n^{\#},G_n)$; $d_{01}(Q_{2n}^{\#},Q_{20})$ and $d_{02}(G_n^{\#},G_0)$ are $O_P(n^{-1/2-\alpha/4})$.

\paragraph{Preservation of rate of TMLE under sampling from $P_n$:}
Lemma \ref{lemmathalmleb} in Appendix  \ref{AppendixC} proves that indeed
$d_{01}(Q_n^{\#*},Q_0)$  converges at same rate $O_P(n^{-1/2-\alpha/4})$ as $d_{01}(Q_n,Q_0)$.

\paragraph{Consistency of nonparametric bootstrap for HAL-TMLE:}
This verifies all conditions of Theorem \ref{thnpboothaltmle} which establishes the asymptotic efficiency and asymptotic consistency of the nonparametric bootstrap.

\begin{theorem}
Consider statistical model ${\cal M}$ defined by (\ref{modelate}), indexed by sectional variation norm bounds $C^u$. Let the  statistical target parameter $\Psi:{\cal M}\rightarrow\openr$ be defined by $\Psi(P)=\Psi_1(P)-\Psi_0(P)$, where $\Psi_a(P)=E_PE_P(Y\mid A=a,W)$.
Consider the HAL-TMLE $Q_n^*$ of $Q_0$ defined above. 
We have that $\Psi(Q_n^*)$ is asymptotically efficient, i.e. $n^{1/2}(\Psi(Q_n^*)-\Psi(Q_0))\Rightarrow_d N(0,\sigma^2_0)$, where
$\sigma^2_0=P_0 \{D^*(P_0)\}^2$.
In addition, conditional on $(P_n:n\geq 1)$, $Z_n^{1,\#}=n^{1/2}(\Psi(Q_n^{\#*})-\Psi(Q_n^*))\Rightarrow_d N(0,\sigma^2_0)$.
This can also be applied to the setting in which $C^u$ is replaced by the cross-validation selector $C_n^u$ defined above, or any other data adaptive selector $\tilde{C}_n^u$ satisfying $P(C_n^u\leq \tilde{C}^u_n\leq C^u)=1$.
\end{theorem}


\subsection{Nonparametric estimation of integral of square of density}
\paragraph{Statistical model, target parameter, canonical gradient:}
Let $O\in \openr^d$ be a multivariate random variable with probability distribution $P_0$ with support $[0,\tau]$.
Let ${\cal M}$ be a  nonparametric model dominated by Lebesgue measure $\mu$, where we assume that for each $P\in {\cal M}$ its density
$p=dP/d\mu$ is bounded from above by some $M<\infty$ and from below by some $\delta>0$. Consider the parametrization $p(x)=p(Q)(x)=c(Q)\frac{1}{1+\exp(-Q(x))}$, where $c(Q)$ is the normalizing constant. Our model ${\cal M}$ also assume that $Q$ varies over all cadlag functions with sectional variation norm bounded by $C^u<\infty$. Due to the $C^u$-bound, we also have that any density $p$ in our model is bounded from above by an $M=M(C^u)$ and from below by a $\delta=\delta(C^u)$. This shows that our model ${\cal M}=\{P: Q(P)\in D[0,\tau], \pl Q\pl_v^*<C^u\}$ is of the type  (\ref{calFmodel}).

An alternative formulation that avoids a normalizing constant $c(Q)$ is the following. For sake of presentation, let's consider the case that $d=2$. We factorize $p(x)=p_1(x_1)p_2(x_2\mid x_1)$. 
Subsequently, we parametrize $p_1$ in terms of its hazard $\lambda_1(x_1)=p_1(x_1)/\int_{x_1}^{\tau_1}p_1(u)du$, and $p_2$ in terms of its conditional hazard $\lambda_2(x_2\mid x_1)=p_2(x_2\mid x_1)/\int_{x_2}^{\tau_2} p_2(v\mid x_1)dv$. We then parametrize $\lambda_1=\exp(Q_1)$ and $\lambda_2=\exp(Q_2)$ so that the functions $Q_1(x_1)$ and $Q_2(x_2\mid x_1)$ are unrestricted. This then defines a parameterization $p=p_{Q_1,Q_2}$. 
The log-likelihood provides valid loss functions for $Q_1$ and $Q_2$, so that the HAL-MLE can be computed by  maximizing the log-likelihood over linear combinations of a large number of indicator basis functions with a bound on the $L_1$-norm of its coefficients. 

The target parameter $\Psi:{\cal M}\rightarrow\openr$ is defined as $\Psi(P)=E_Pp(O)=\int p^2(o)d\mu(o)$. We can represent $\Psi(P)$ as a function of $Q$ or the density $p$, so that we will also denote it with $\Psi(Q)$ or $\Psi(p)$.
This target parameter is pathwise differentiable at $P$ with canonical gradient 
\[
D^*(Q)(O)=2 (p(O)-\Psi(p)),\]
where $p=p(Q)$. 

\paragraph{Exact second order remainder:}
It implies the following exact second-order expansion:
\[
\Psi(Q)-\Psi(Q_0)=(P-P_0)D^*(Q)+R_{20}(Q,Q_0),\]
where 
\[
R_{20}(Q,Q_0)\equiv -\int (p-p_0)^2 d\mu .\]
\paragraph{Loss function:}
As loss function for $Q$ we could consider the log-likelihood loss $L(Q)(O)=-\log p(O)$ with $d_0(Q,Q_0)=P_0\log p_0/p$, where again $p_0=p(Q_0)$ and $p=p(Q)$.
We have $\pl p^{1/2}-p_0^{1/2}\pl_{P_0}^2\leq  P_0 \log p_0/p $ so that 
\begin{eqnarray*}
| R_{20}(Q,Q_0)| &=&\int (p-p_0)^2 d\mu \\
&=&\sup\frac{(p^{1/2}+p_0^{1/2})^{2}}{p_0} \int (p^{1/2}-p_0^{1/2})^2 dP_0\\
&\leq&  M/\delta P_0\log p_0/p =M/\delta d_0(Q,Q_0).
\end{eqnarray*}
Alternatively, we could consider the loss function 
\[
L(Q)(O)=-2p(O)+\int p^2d\mu .\]
Note that this is  indeed a valid loss function with loss-based dissimilarity given by \begin{eqnarray*}
d_0(Q,Q_0)&=&P_0 L(Q)-P_0L(Q_0)\\
&=&-2\int p(o)p_0(o)d\mu(o)+\int p^2d\mu+2\int p_0^2 d\mu-\int p^2_0 d\mu\\
&=& \int (p-p_0)^2 d\mu .\end{eqnarray*}
\paragraph{Bounding second order remainder:}
Thus, if we select this loss function, then we have
\[
| R_{20}(Q,Q_0)| =d_0(p,p_0) .\]
In terms of our general notation, we now have $f(x)=x^2$ for the upper bound on $R_{20}$ so that $| R_{20}(Q,Q_0)| =f(d_0^{1/2}(Q,Q_0))$.
The canonical gradient is indeed continuous in $Q$ as stated in (\ref{contDstar}) and the bounds $M_1,M_2,M_3$ (\ref{sectionalvarbound}) are obviously finite and can be expressed in terms of $(C^u,M(C^u),\delta(C^u))$.
This verifies the assumptions on our model as stated in Section 2.

\paragraph{HAL-MLE and CV-HAL-MLE:}
Let $Q_n=\argmin_{Q,\pl Q\pl_v^*<C^u}P_n L(Q)$ be the HAL-MLE.
where $Q$ can be represented by our general representation (\ref{Frepresentation}), $Q(o)=Q(0)+\sum_{s\subset \{1,\ldots,d\}}\int_{(0_s,o_s]} dQ_s(u_s)$,
and constrained to satisfy 
\begin{align*}    
| Q(0)| + \sum_{s\subset\{1,\ldots,d\}}\int_{(0_s,\tau_s]} |dQ_s(u_s)| \leq C^u.
\end{align*}
Let's denote this $Q_n$ with $Q_{n,C^u}$.
Thus, for a given $C$, computation of $Q_{n,C}$ can be done with a LASSO type algorithm. Let $C_n=\argmin_CE_{B_n}P_{n,B_n}^1L(\hat{Q}_C(P_{n,B_n}^0))$ be the cross-validation selector of $C$, as defined in previous example. If we set $C=C_n$, then we obtain the CV-HAL-MLE $Q_n=Q_{n,C_n}$.
By our general result on HAL-MLE for bounded loss functions, we have (for both the log-likelihood loss and $L^2$-loss functions) $d_0(Q_n,Q_0)=O_P(n^{-1/2-\alpha/4})$.

\paragraph{TMLE using Local least favorable submodel and log-likelihood loss:}
Let $p_n=p(Q_n)$. A possible local least favorable submodel through $p_n$ when using the log-likelihood loss is given by  $p_{n,\epsilon}^{lfm}=(1+\epsilon D^*(p_n)) p_n$ for $\epsilon$ in a small enough neighborhood so that $p_{n,\epsilon}>0$ everywhere: for example, $|\epsilon|<1/\pl p_n\pl_{\infty}$.
Let $\epsilon_n=\argmin_{\epsilon}P_n L(p_{n,\epsilon}^{lfm})$. If $\epsilon_n$ is not in the interior, then one would set $p_n^1=p_{n,\epsilon_n}$ and iterate this updating process till $\epsilon_n$ falls in interior. The final update is denoted with $p_n^*$. As sample size increases, with probability tending to one $\epsilon_n$ will already be in the interior, so that $p_n^*$ would be a closed form one-step TMLE.  Due to the $o_P(n^{-1/4})$-rate of convergence of the HAL-MLE $p_n$,  it follows that $P_n D^*(p_{n,\epsilon}^{lfm})=o_P(n^{-1/2})$. 

\paragraph{TMLE using Universal least favorable submodel and log-likelihood loss:}
One can also define  a universal least favorable submodel \citep{ULFM16} by recursively applying the above local least favorable submodel:
\[
p_{n,\epsilon+d\epsilon}=p_{n,\epsilon,d\epsilon}^{lfm},\]
where $p_{n,\epsilon,d\epsilon}^{lfm}$ is the local least favorable submodel through $p_{n,\epsilon}$ at parameter value $d\epsilon$. In this manner, the local moves of the local least favorable submodel describe a submodel satisfying $\frac{d}{d\epsilon}\log p_{n,\epsilon} =D^*(p_{n,\epsilon})$ at each $\epsilon$. Again, let $\epsilon_n=\argmin_{\epsilon}P_n L(p_{n,\epsilon})$ and $p_n^*=p_{n,\epsilon_n}$, which now satisfies $P_n D^*(p_n^*)=0$ exact.

\paragraph{HAL-TMLE:}
The TMLE of $\Psi(Q_0)=\Psi(p_0)$ is the plug-in estimator $\psi_n^*=\Psi(p_n^*)=\int p_n^{*2}d\mu$.
Lemma \ref{lemmathalmle} in  Appendix \ref{AppendixA} proves $d_{01}(Q_n^*,Q_0)=O_P(n^{-1/2-\alpha/4})$.

\paragraph{Efficiency of HAL-TMLE and CV-HAL-TMLE:}
Theorem \ref{thefftmle} shows that $\Psi(p_n^*)$ is asymptotically efficient, where one can either choose the HAL-MLE with fixed index $C=C^u$  or one can set $C=C_n$ equal to cross-validation selector defined above (or any other data adaptive selector $\tilde{C}_n$ with $P(C_n\leq \tilde{C}_n\leq C^u)=1$.

\paragraph{Asymptotic validity of the nonparametric bootstrap for the HAL-MLE:}
Let $C$ be given. 
As remarked in the previous example, computation of the HAL-MLE  $Q_n^{\#}=\argmin_{\pl Q\pl_v^*\leq C,Q\ll Q_n}P_n^{\#}L(Q)$ is much faster than the computation of $Q_n=\argmin_{\pl Q\pl_v^*\leq C}P_nL(Q)$, due to only having to minimize the empirical risk over the bootstrap sample over the linear combinations of indicator functions that had non-zero coefficients in $Q_n$. 
By Theorem \ref{thnpbootmle} it follows that $d_{n}(Q_n^{\#},Q_n)$ and $d_0(Q_n^{\#},Q_0)$ are $O_P(n^{-1/2-\alpha/4})$.
For example, if we use the $L^2$-loss function above, then $d_0(Q,Q_0)=\int (p-p_0)^2 d\mu$,
and, we note also that
\begin{eqnarray*}
d_n(Q_n^{\#},Q_n)&=&P_n\{ L(Q_n^{\#})-L(Q_n)\}\\
&=&P_n\{ -2(p_n^{\#}-p_n)+\int p_n^{\#2} d\mu-\int p_n^{2} d\mu \}\\
&=&P_n\{ -2(p_n^{\#}-p_n)+\int (p_n^{\#}-p_n)(p_n^{\#}+p_n) d\mu \} \\
&=&P_n\{ \int(-2+2p_n) (p_n^{\#}-p_n)d\mu  \}\\
&=&\int (p_n^{\#}-p_n)^2 d\mu  .
\end{eqnarray*}
This shows that the empirical dissimilarity also equals the square of an  $L^2$-norm.
Thus, application of Theorem \ref{thnpbootmle} now shows that $\int (p_n^{\#}-p_n)^2 d\mu$  and $\int (p_n^{\#}-p_0)^2 d\mu$ are both $O_P(n^{-1/2-\alpha/4})$.

\paragraph{Preservation of rate for HAL-TMLE under sampling from $P_n$:}
Lemma \ref{lemmathalmleb} in  Appendix \ref{AppendixC} establishes that $d_{01}(Q_n^{\#*},Q_0)=O_P(n^{-1/2-\alpha/4})$.

\paragraph{Asymptotic consistency of the bootstrap for the HAL-TMLE:}
This verifies all conditions of Theorem \ref{thnpboothaltmle} which establishes the asymptotic efficiency and asymptotic consistency of the nonparametric bootstrap.

\begin{theorem}
Consider the model ${\cal M}$ defined by upper  bound $C^u<\infty$ on the sectional variation norm of $Q$ over $[0,\tau]$. Let $\Psi(Q)=\int p(Q)^2d\mu$, which is also denote with $\Psi(p)$.
Consider the one-step TMLE based on the local least favorable submodel or universal least favorable submodel and the log-likelihood loss. 

We have that $\Psi(p_n^*)$ is asymptotically efficient, i.e. $n^{1/2}(\Psi(p_n^*)-\Psi(p_0))\Rightarrow_d N(0,\sigma^2_0)$, where
$\sigma^2_0=P_0 \{D^*(P_0)\}^2$.

In addition, conditional on $(P_n:n\geq 1)$, $Z_n^{1,\#}=n^{1/2}(\Psi(p_n^{\#*})-\Psi(p_n^*))\Rightarrow_d N(0,\sigma^2_0)$.

This theorem can also be applied to the setting in which $C^u=C_n$.
\end{theorem}

\section{Simulation study evaluating performance of bootstrap method}\label{sectsix}

\subsection{Average treatment effect}

To illustrate the finite sample performance of the proposed bootstrap method, we
simulate a continuous outcome $Y$, a binary treatment $A$, and a continuous
covariate $W$ that confounds $Y$ and $A$. The random variables are drawn from a
family of distributions indexed by $a_1$, which characterizes the conditional
distribution of $Y$, given $A$ and $W$. The distribution of variables are as
follows: $W \sim N(0, 4^2, -10, 10)$ is drawn i.i.d. from a truncated normal
distribution with mean equals 0, standard deviation 4, bounded within
$[-10,10]$. $A \sim Bernoulli(\bar{G}(W))$ is a Bernoulli binary random variable, with
a probability $\bar{G}(W)$ as a function of $W$, given by
\begin{align*}
    \bar{G}(W)= 0.3 + \min(0.1 W \sin(0.1 W) + \epsilon_1, 0.4) 
\end{align*}
where $\varepsilon_1 \sim N(0, 0.05^2)$. $Y = 3\sin(a_1W) + A + \varepsilon_2$ is
a sinusoidal function of $W$, where $\varepsilon_2 \sim N(0, 1)$, which defines $\bar{Q}_0(A,W)=3\sin(a_1 W)+A$. $a_1$ controls
the amplitude of the sinusoidal function. Increasing $a_1$
(frequency) of the sin function increases the sectional variation norm of $\bar{Q}_0$ proportionally, so that
estimating $Q_0$ becomes more difficult under fixed sample size. In our study,
we increase $a_1$ while fixing the sample size and fixing the HAL-MLE $G_n$ of $G_0$, so that the second order
remainder increases in magnitude.
The value of the parameter of interest, ATE $\psi_0=\Psi(P_0)$, is
1. The experiment is repeated 1000 times. The estimation routine
including the tuning parameter search is implemented in the \texttt{ateBootstrap}
function in the open-source R package \texttt{TMLEbootstrap} \citep{cai2018software}. A thousand replications of the simulation 1 under sample size 100 and 200 bootstrap repetitions takes 12 CPU hour on an Intel Core i7 4980HQ CPU.

\begin{figure}[!ht]
    \centering
    \includegraphics[width=0.9\textwidth]{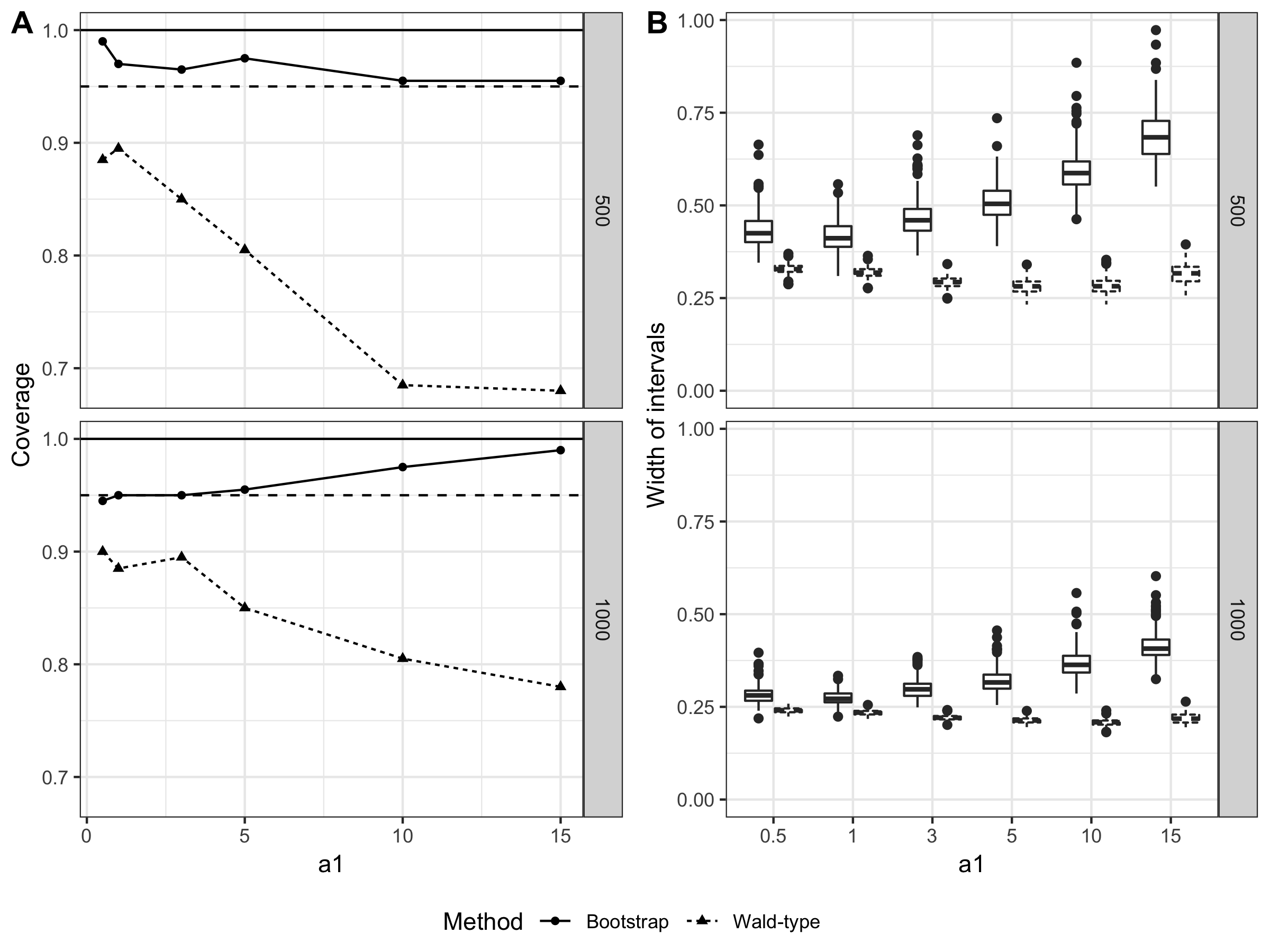}
    \caption{Results for ATE parameter comparing our bootstrap method and
    classic Wald-type method as a function of the $a_1$ coefficient (sectional
    variation norm) of the $\bar{Q}_0$ function. Panel A is the coverage of the
    intervals, where dashed line indicate 95\% nominal coverage. Panel B is the
    widths of the intervals. Within each panel, the upper plot is under sample
    size 500 and the lower plot is under sample size 1000.}
    \label{fig:ate-simu}
\end{figure}

To analyze the above simulated data, we compute the coverage and width of confidence interval  of the Wald-type confidence interval where the nuisance functions $(\bar{Q}_0, G_0)$ are estimated using HAL-MLE($\lambda_{CV}$) and nonparametric bootstrap confidence interval presented  in Section \ref{sectfour}, where the choice $\lambda$ in $\bar{Q}_{n,\lambda}^*$ is set equal to the plateau selector $\lambda_{plateau}$.  Recall that the nonparametric bootstrap of the HAL-TMLE at this choice $\lambda_{plateau}$ is used to determine the quantiles for the confidence interval around the TMLE $\Psi(Q_n^*)$. Wald-type interval reflects common practice for statistical inference based on the TMLE.
Results under samples sizes 500 and 1000 are shown in Figure \ref{fig:ate-simu}.

\begin{figure}[!ht]
    \centering
    \includegraphics[width=1\textwidth]{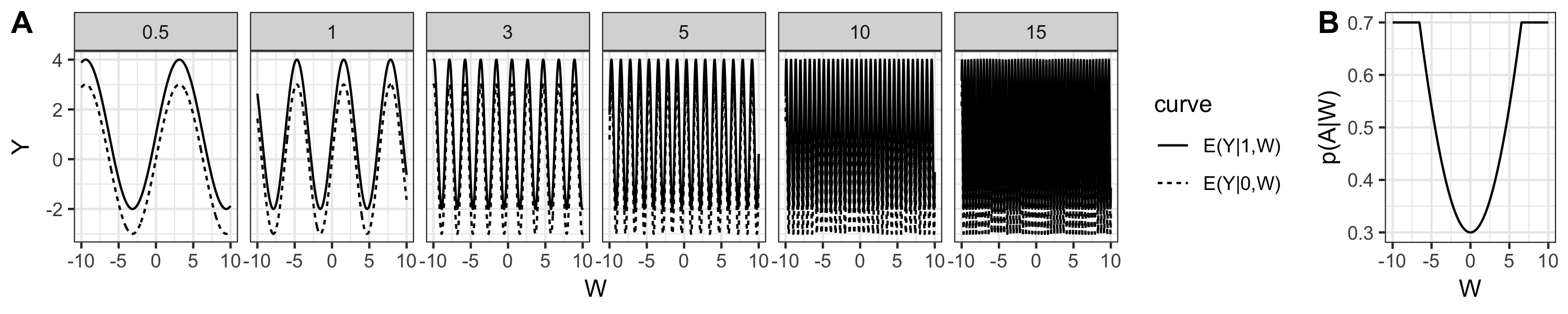}
    \caption{(A) True conditional expectation functions of outcome $E(Y|A=1,W)$ and $E(Y|A=0,W)$ at $a_1 = 0.5,1,3,5,10,15$ and (B) true propensity score function}
    \label{fig:ate-dgd}
\end{figure}

The simulation results reflect what is expected based on theory. In particular,
as the sectional variation norm of the $\bar{Q}_0$ becomes large relative to
sample size, the HAL regression fit of $\bar{Q}_0$ in the finite sample is not ideal, which leads to low coverage of Wald-type interval. On the other hand, the bootstrap confidence  intervals reflect the deteriorating  second-order remainder in the sampling distribution of the HAL-TMLE of $\bar{Q}_0$, and, as a result, the coverage is very close to nominal and is
robust to increasing sectional variation norm ($a_1$). The results for sample
size 1000 confirm our asymptotic analysis of the methods, with Wald-type coverage
improving and two methods eventually converging to nominal coverage.

\subsection{Average density value}

\begin{figure}[!ht]
    \centering
    \includegraphics[width=1\textwidth]{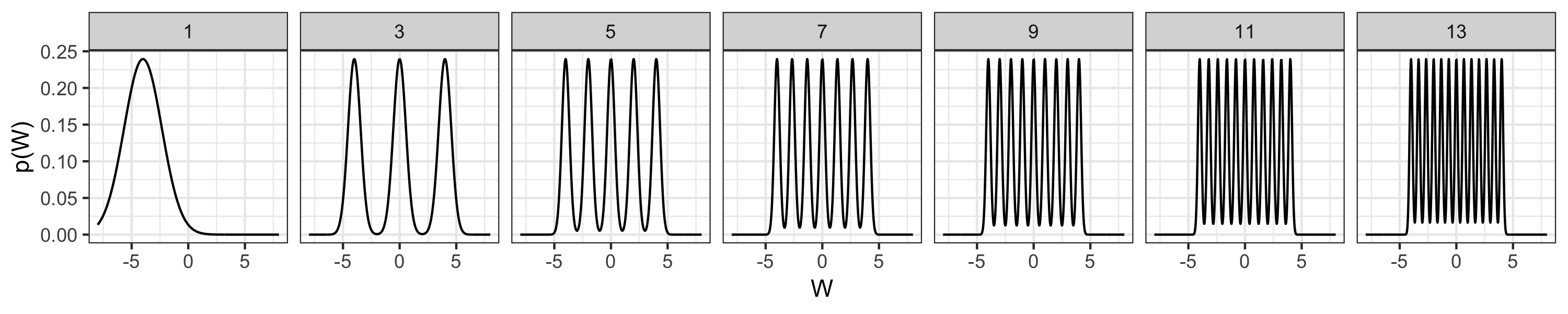}
    \caption{True probability density function $f(x;\theta_K)$ at $K = 1,3,5,7,9,11,13$}
    \label{fig:avgdens-dgd}
\end{figure}

As we demonstrated, this problem has a non-forgiving  second-order
remainder term that is proportional to the $L^2$-norm of $p_n^*-p_0$, which makes this example very suitable for evaluating finite sample coverage of the bootstrap methods.  To illustrate our proposed method and
explore finite-sample performance, we simulate a family of univariate densities
with increasing sectional variation norm. $$f(x;\theta_K) =
\frac{1}{K}\sum\limits_{k = 1}^K{g(x;\mu_k,\sigma_K)},$$ where
$$g(x;\mu_k,\sigma_K) = \frac{1}{\sqrt{2\pi}
\sigma_K}\exp[-\frac{1}{2}(x-\mu_k)^2/\sigma_K^2].$$ For a given $K$, $\mu_k, k
= 1,...,K$ are equi-distantly placed in interval $[-4,4]$. $\sigma_K=10/K/6$.
The true sectional variation norm of the density increases roughly linearly with
K, that is $\|f_K\|_v^* = K\|f_1\|_v^*, K = 1,...,13$. Examples of the density family
for $K$ values used in the simulation are shown in Figure \ref{fig:avgdens-dgd}.
We simulate from univariate densities for the sake of presentation and we expect
our results to be informative for higher dimensional densities as well, since the main difference will be that a  sectional
variation norm of a multivariate function is generally larger than that of a  univariate function. 
The value of the parameter of interest $\Psi(p_0)=\int p_0^2dx$ does not change much as a function of the choice  $K$ of data distribution. For each data distribution,  the
experiment is repeated 1000 times. As in our ATE simulation, we compute
coverages and widths  of the Wald-type confidence intervals (that ignore second order remainders) and our 
HAL-TMLE bootstrap confidence intervals using our plateau selector $\lambda_{plateau}$.

We parametrized the density in terms of its hazard, discretized the hazard making it piecewise constant across a large number of bins (like histogram density estimation), parametrized this piecewise constant hazard with a logistic regression for the probability of falling in bin $h$, given it exceeded bin $h-1$. We fitted this hazard with a logistic regression based HAL-MLE using the longitudinal data format common for hazard estimation  (i.e., an observation $O_i$ is coded by a number of rows with binary outcome equal to zero and a final row with outcome $1$). 
The HAL-MLE of this hazard yields the corresponding HAL-MLE of the density itself.
The HAL-TMLE updates the HAL-MLE density estimator with a TMLE update using the universal least favorable submodel and log-likelihood loss. 
The software implementations can be found in the \texttt{cv\_densityHAL}
function in the open-source R package \texttt{TMLEbootstrap} \citep{cai2018software}.

\begin{figure}[!ht]
    \centering
    \includegraphics[width=0.9\textwidth]{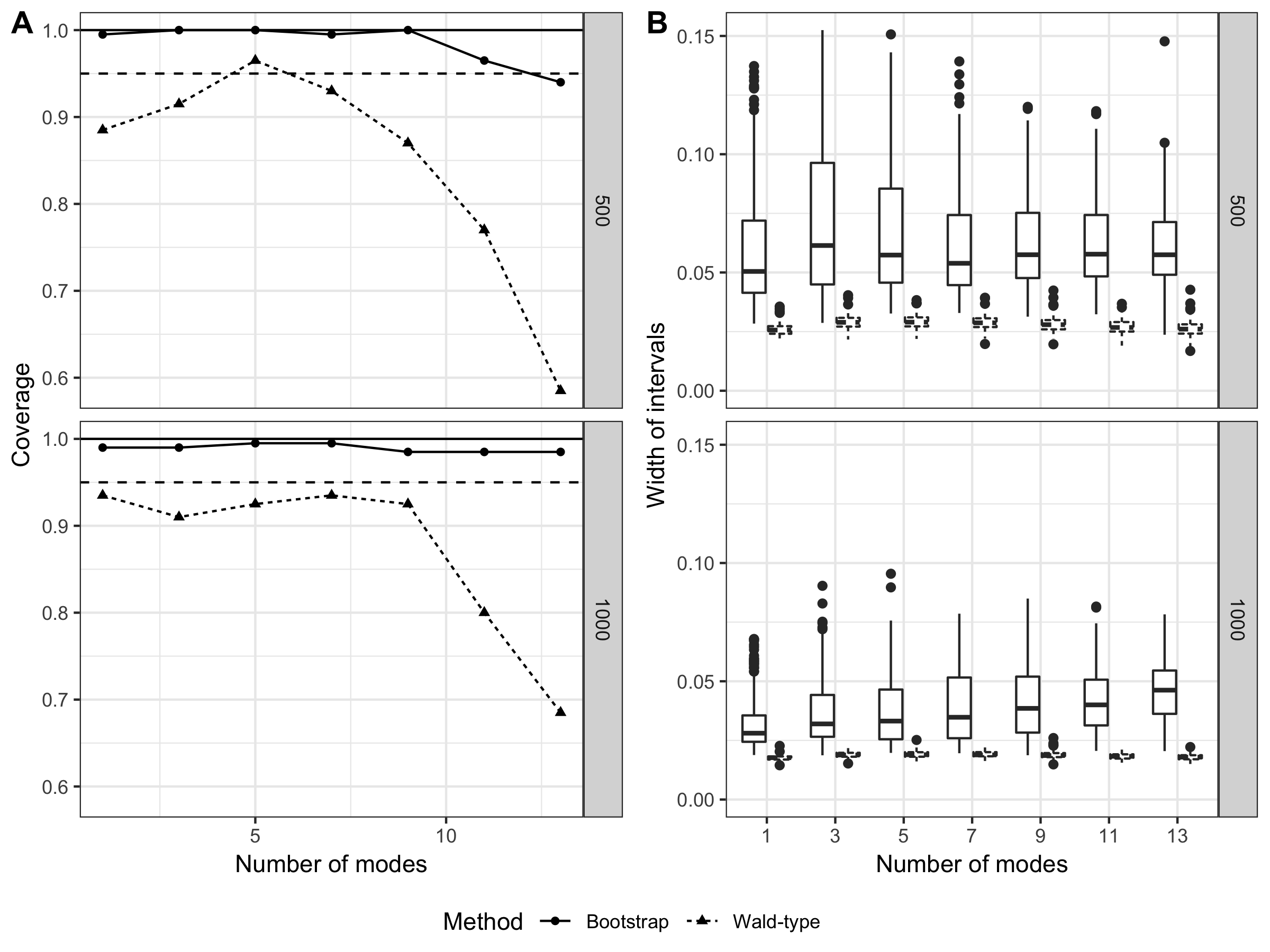}
    \caption{Results for average density value parameter comparing our bootstrap
    method and classic Wald-type method as a function of the number of modes in
    true density (sectional variation norm). Panel A is the coverage of the
    intervals, where dashed line indicate 95\% nominal coverage. Panel B is the
    widths of the intervals. Within each panel, the upper plot is under sample
    size 500 and the lower plot is under sample size 1000.}
\end{figure}

The simulations reflect what is expected based on theory: the bootstrap confidence interval
has superior coverage relative to the  Wald-type confidence interval, uniformly across
different sample sizes and data distributions. In particular, as the true sectional variation norm
increases (with the number of modes in the density), the  second-order
remainder term increases so that the  Wald-type interval coverage declines. On the other
hand, the bootstrap confidence intervals  reflect the behavior of the second order remainder and thereby  increase in width as the performance of the HAL-MLE deteriorates (due to increased complexity of true density). The bootstrap
confidence interval controls the coverage close to the nominal rate and its coverage is not very 
sensitive to the true sectional variation norm of the density function. When
sample size increases to 1000, the Wald-type interval coverage increases, and in
simple cases where the true sectional variation norm is small, Wald-type
coverage reaches its desired nominal covarage. 

\section{Discussion}\label{sectdisc}
On one hand, in parametric models and, more generally, in  models small enough so that  the MLE is still well behaved, one can use the nonparametric bootstrap to estimate the sampling distribution of the MLE. It is generally understood that in these small models the nonparametric bootstrap  outperforms estimating the sampling distribution with a normal distribution (e.g., with variance estimated as the sample variance of the influence curve of the MLE), by picking up the higher order behavior of the MLE, {\em if asymptotics has not set in yet}. In such small models, reasonable sample sizes already achieve the normal approximation in which case the Wald type confidence intervals will perform well. 
Generally speaking, the nonparametric bootstrap is a valid method when the estimator is a compactly differentiable function of the empirical measure, such as the Kaplan-Meier estimator (i.e., one can apply the functional delta-method to analyze such estimators) \citep{Gill89}\citep[Theorem 3.9.11 in][]{vanderVaart&Wellner96}. These are estimators that essentially do not use  smoothing of any sort.

On the other hand, efficient estimation of a pathwise differentiable target parameter in  large realistic models generally requires estimation of the data density, and thereby machine learning such as super-learning to estimate the relevant parts of the data distribution. Therefore, efficient one-step estimators or TMLEs are not compactly differentiable functions of the data distribution.
Due to this reason, we moved away from using the nonparametric bootstrap to estimate its sampling distribution, since it represents a generally inconsistent method (e.g., a cross-validation selector behaves very differently under sampling from the empirical distribution than under sampling from the true data distribution) \citep{coyle2018targetboot}. Instead we estimated the normal limit distribution by estimating the variance of the influence curve of the estimator.  

Such an influence curve based method is asymptotically consistent and therefore results in asymptotically valid $0.95$-level confidence intervals. However,
in such large models the nuisance parameter estimators will converge at slow rates (like $n^{-1/4}$ or slower) with large constants depending on the size of the model, so that for normal sample sizes the exact second-order remainder could be easily larger than the leading empirical process term with its normal limit distribution. 
So one has to pay a significant price for using the computationally attractive influence curve based confidence intervals, where inference is ignoring the remainder terms $(P_n-P_0)(D^*_n - D^*_0)$ and $R_2(P_n^*,P_0)$. In finite sample these remainder terms can have non-zero expectation or have a large variance, so the influence curve-based inference using a normal limit distribution can be off-centered or less spread out than the actual sampling distribution of the estimator.


One might argue that one should use a model based bootstrap instead by sampling from an  estimator of the density of the data distribution. General results show that such a  model based bootstrap method will be asymptotically valid as long as the density estimator is consistent \citep{arcones1989boot_mean,gine1989necessary_boot_mean,arcones1992boot_mestimator}. This is like carrying out a simulation study for the estimator in question using an estimator of the true data distribution as sampling distribution. However, estimation of the actual density of the data distribution is itself a very hard problem, with bias heavily affected by the curse of dimensionality, and, in addition, it can be immensely burdensome to construct such a density estimator and sample from it when the data is complex and high dimensional.

As demonstrated in this article, the HAL-MLE provides a solution to this bottleneck. The HAL-MLE($C^u$)  of the nuisance parameter is an actual MLE minimizing the empirical risk over a infinite dimensional parameter space (depending on the model ${\cal M}$) in which it is assumed that the sectional variation norm of the nuisance parameter is bounded by universal constant $C^u$. This MLE is still well behaved by being consistent at a rate that is in the worst case still faster than $n^{-1/4}$. However, this MLE is not an interior MLE, but will be on the edge of its parameter space: the MLE will itself have sectional variation norm equal to the maximal allowed value $C^u$. Nonetheless, our analysis shows that it is still  a smooth enough function of the data (while not being compactly differentiable at all) that it is equally well behaved under sampling from the empirical distribution.  

As a consequence of this robust behavior of the HAL-MLE, for  models in which the nuisance parameters of interest are cadlag functions with a universally bounded sectional variation norm (beyond possible other assumptions), we presented asymptotically consistent estimators of the sampling distribution of the HAL-TMLE  of the target parameter of interest using the nonparametric bootstrap. 

Our estimators of the sampling distribution are highly sensitive to the curse of dimensionality, just as the sampling distribution of the HAL-TMLE itself: specifically, the HAL-MLE on a bootstrap sample will  converge just as slowly to its truth as under sampling from the true distribution. Therefore, in high dimensional estimation problems, we  expect highly significant gains in valid inference relative to Wald type confidence intervals  that are purely  based on the normal limit distribution of the HAL-TMLE.

In general, the  user will typically not know how to select the upper bound $C^u$ on the sectional variation norm of the nuisance parameters (except if the nuisance  parameters are cumulative distribution functions). Therefore, for the sake of estimation of $Q_0$ and $G_0$ we recommend to select this bound with cross-validation. Due to the oracle inequality for the cross-validation selector $C_n$ (which only relies on a bound on the supremum norm of the loss function), the data adaptively selected upper bound will be selected larger than (but close to) the true sectional variation norm  $C_0$ of the nuisance parameters $(Q_0,G_0)$, as sample size increases. 

Even though, for this cross-validation selector $C_n$,  our bootstrap estimators will still be guaranteed to be consistent for its normal limit distribution,  this choice $C_n$ will be trading off bias and variance for the sake of estimation of the nuisance parameter. As a consequence, in practice this $C_n^u$ might often end up selecting a value significantly smaller than the true sectional variation norms of $Q_0$ and $G_0$. This is comparable with selecting a models for $Q_0$ and $G_0$ to be used in the bootstrap that are potentially much smaller than a model that would be needed to capture $Q_0$ and $G_0$. That is, our proposed bootstrap would then still not capture the full complexity of the estimation problem and still result in anti-conservative confidence intervals. 
Therefore we proposed a finite sample modification for the nonparametric bootstrap of the HAL-TMLE by using the bootstrap distribution for the HAL-TMLE at fixed sectional variation norm determined by a plateau selector (instead of the cross-validation selector of $C^u$). Our proposed finite sample modification also uses a scaling $\sigma_n^2$ that incorporates both  bias and variance of the bootstrap distribution. Our simulations demonstrate the importance of this finite sample modification and showcases excellent finite sample coverage. Any improvements for variance estimation relative to using the empirical variance of the influence curve can be incorporated naturally in this method, such as  the plug-in variance estimator that is robust under data sparsity presented in \citep{Tranetal18}. 

There are a number of important future directions to this research. One direction is to derive finite-sample bounds on our bootstrap interval coverage probability, which will give additional guarantees for applications.

\subsection*{Acknowledgement.}
This research is funded by NIH-grant 
5R01AI074345-07.

\bibliography{combined,TLB2}


\appendix 
\section*{Appendix.}

The HAL-MLEs on the original sample and bootstrap sample will be defined below as
$Q_n=\arg\min_{Q\in Q({\cal M})}P_n L_1(Q)$ and $Q_n^{\#}=\arg\min_{Q\in Q({\cal M})}P_n^{\#}L_1(Q)$, and, if we assume the extra structure (\ref{calFmodel}) so that we know that $Q({\cal M})$ is itself defined as a space of cadlag functions with bounds on its sectional variation norm, then we let 
$Q_n^{\#}=\arg\min_{Q\in Q({\cal M}), Q\ll Q_n,\pl Q\pl_v^*\leq \pl Q_n\pl_v^*}P_n^{\#}L_1(Q)$.
In general, one can add restrictions to the parameter space over which one minimizes in the definition of $Q_n$ and $Q_n^{\#}$ as long as one guarantees that, with probability tending to 1, $P_nL_1(Q_n)\leq P_nL_1(Q_0)$, and, with probability tending to 1, conditional on $P_n$,  $P_n^{\#}L_1(Q_n^{\#})\leq P_n^{\#}L_1(Q_n)$. For example, this allows one to use an upper bound $C_n^u$ on the sectional variation norm  in the definition of $Q_n$ if we know that $C_n^u$ will be larger than the true $C_0^u=\pl Q_0\pl_v^*$ with probability tending to 1.

\section{Proof that the one-step TMLE $Q_n^*$ preserves rate of convergence of $Q_n$}\label{AppendixA}
The following lemma establishes that the one-step TMLE $Q_n^*=Q_{n,\epsilon_n}$ preserves the rate of convergence of $Q_n$, where $Q_{\epsilon}$ is a univariate local least favorable submodel through $Q$ at $\epsilon =0$.
We use the  notation $L_1(Q_1,Q_2)=L_1(Q_1)-L_1(Q_2)$.
\begin{lemma}\label{lemmathalmle}
Let $Q_n=\arg\min_{Q\in Q({\cal M})}P_n L_1(Q)$ be the HAL-MLE of $Q_0$,  and let $\epsilon_n =\arg\min_{\epsilon}P_n L_1(Q_{n,\epsilon})$ for some parametric (e.g, local least favorable) submodel $\{Q_{\epsilon}:\epsilon\}\subset {\cal M}$.
 Assume the bounds (\ref{M20bounds})$, (\ref{sectionalvarbound})$ on loss function $L_1(Q)$, so that we also know $d_{01}(Q_n,Q_0)=O_P(n^{-1/2-\alpha/4})$. 

Then, 
\begin{equation}\label{thalmle}
d_{01}(Q_n^*,Q_0)=O_P(n^{-1/2-\alpha/4}).\end{equation}
Specifically, we have 
\[
d_{01}(Q_n^*,Q_0) \leq -(P_n-P_0)L_1(Q_{n,\epsilon_{n}},Q_n)+d_{01}(Q_n,Q_0).\]
\end{lemma}

This also proves that the $K$-th step TMLE using a {\em finite} $K$ (uniform in $n$) number of iterations satisfies $d_{01}(Q_n^*,Q_0)\leq d_{01}(Q_n,Q_0)+O_P(n^{-1/2-\alpha/4})$. So  if $r_n=P_n D^*(Q_{n,\epsilon_n},G_n)$ is not yet $o_P(n^{-1/2})$, then  one should consider  a $K$-th step  TMLE to guarantee that $r_n$ is small enough to be neglected (we know that the fully iterated TMLE will solve $P_n D^*(Q_n^*,G_n)=0$, but this one is harder to analyze).
\newline
{\bf Proof of Lemma \ref{lemmathalmle}:}
 We have
\begin{eqnarray*}
P_0L_1(Q_n^*)-P_0L_1(Q_0)
&=& P_0L_1(Q_{n,\epsilon_{n}},Q_n) +P_0L_1(Q_n,Q_0)\\
&=&(P_0-P_n)L_1(Q_{n,\epsilon_{n}},Q_n) +P_n L_1(Q_{n,\epsilon_{n}},Q_n)\\
&&+d_{01}(Q_n,Q_0)\\
&\leq &-(P_n-P_0)L_1(Q_{n,\epsilon_{n}},Q_n)+d_{01}(Q_n,Q_0)\\
\end{eqnarray*}
Since $L_1(Q_{n,\epsilon_n},Q_n)$ falls in class of cadlag functions with a universal bound on sectional variation norm (i.e., a Donsker class), and $d_{01}(Q_n,Q_0)=O_P(n^{-1/2-\alpha/4})$, it follows that $d_{01}(Q_n^*,Q_0)=O_P(n^{-1/2})+O_P(n^{-1/2-\alpha/4})$.
Now, we use that the $L^2(P_0)$-norm of $L_1(Q_{n,\epsilon_n},Q_n)$ is bounded by sum of $L^2(P_0)$-norm of $L_1(Q_n^*)-L_1(Q_0)$ and $L_1(Q_n)-L_1(Q_0)$. These latter $L^2(P_0)$-norms can be bounded by $d_{01}^{1/2}(Q_n^*,Q_0)$ and $d_{01}^{1/2}(Q_n,Q_0)$, which thus converges at rate $O_P(n^{-1/4-\alpha/2})$.  Again, by empirical process theory, using that we now know $P_0\{ L_1(Q_n^*,Q_n)\}^2=o_P(n^{-1/4-\alpha/2})$, it follows immediately that $d_{01}(Q_n^*,Q_0)=o_P(n^{-1/2})$, but, by using the actual rate for this $L^2(P_0)$-norm, as in the Appendix in \citep{vanderLaan15} it follows that $d_{01}(Q_n^*,Q_0)=O_P(n^{-1/2-\alpha/4})$. 
$\Box$

\section{Asymptotic convergence of bootstrapped HAL-MLE: Proof of Theorem \ref{thnpbootmle}.}\label{AppendixB}

 Theorem \ref{ThB} below shows that both $d_{n1}(Q_n^{\#},Q_n)=P_n \{L_1(Q_n^{\#})-L_1(Q_n)\}$ and $d_{01}(Q_n^{\#},Q_0)$ converge at rate $n^{-1/2-\alpha/4}$. The analogue results apply to $G_n^{\#}$.

\begin{theorem}\label{ThB}
Consider a statistical model ${\cal M}$ 
satisfying (\ref{sectionalvarbound}), (\ref{M20bounds}) on $L_1(Q)$.
Let $Q_n=\arg\min_{Q\in Q({\cal M})}P_nL_1(Q)$ and $Q_n^{\#}=\arg\min_{Q\in Q({\cal M})}P_n^{\#}L_1(Q)$. In a model with extra structure (\ref{calFmodel}) we define \[
Q_n^{\#}=\arg\min_{Q\in Q({\cal M}),Q\ll Q_n, \pl Q_n^{\#}\pl_v^*\leq \pl Q_n\pl_v^*}P_n^{\#}L_1(Q).\] 
Then, \[ d_{n1}(Q_n^{\#},Q_n)=O_P(n^{-1/2-\alpha/4})\mbox{ and } d_{01}(Q_n^{\#},Q_0)=O_P(n^{-1/2-\alpha/4}).
\]
\end{theorem}
{\bf Proof of Theorem \ref{ThB}:}
We have
\begin{eqnarray}
0&\leq & d_{n1}(Q_n^{\#},Q_n)\equiv P_n\{ L_1(Q_n^{\#})-L_1(Q_n)\}\nonumber \\
&=& -(P_n^{\#}-P_n)\{L_1(Q_n^{\#})-L_1(Q_n)\}+P_n^{\#}\{L_1(Q_n^{\#})-L_1(Q_n)\}\nonumber \\
&\leq& -(P_n^{\#}-P_n)\{L_1(Q_n^{\#})-L_1(Q_n)\}.\label{boota}
\end{eqnarray}
As a consequence, by empirical process theory \citep{vanderVaart&Wellner11}, we have $d_{n1}(Q_n^{\#},Q_n)=P_n L_1(Q_n^{\#})-P_nL_1(Q_n)$ is $O_P(n^{-1/2})$. 
We now note that
\begin{eqnarray}
d_{01}(Q_n^{\#},Q_0)&=&P_0L_1(Q_n^{\#},Q_n)+P_0L_1(Q_n,Q_0)\nonumber\\
&=&(P_0-P_n)L_1(Q_n^{\#},Q_n)+P_n L_1(Q_n^{\#},Q_n)+d_{01}(Q_n,Q_0)\nonumber \\
&=&(P_0-P_n)L_1(Q_n^{\#},Q_n)+d_{n1}(Q_n^{\#},Q_n)+d_{01}(Q_n,Q_0).\label{handy}
\end{eqnarray}
Thus, it also follows that $d_{01}(Q_n^{\#},Q_0)=O_P(n^{-1/2})$.
By assumption \ref{M20bounds} this implies
$P_0\{L_1(Q_n^{\#})-L_1(Q_0)\}^2=O_P(n^{-1/2})$. Note now that $L_1(Q_n^{\#}) - L_1(Q_n) = L_1(Q_n^{\#}) - L_1(Q_0) + L_1(Q_0) - L_1(Q_n)$, using that $P_0\{L_1(Q_n) - L_1(Q_0)\}^2=O_P(n^{-1/2})$, it follows that also $P_0\{L_1(Q_n^{\#})-L_1(Q_n)\}^2=O_P(n^{-1/2})$. By Lemma \ref{lemmadndo}, it follows that also $P_n\{L_1(Q_n^{\#})-L_1(Q_n)\}^2=O_P(n^{-1/2})$. 
With this result in hand, using \citep{vanderVaart&Wellner11} as in Appendix in \citep{vanderLaan15}, it follows that $-(P_n^{\#}-P_n)\{L_1(Q_n^{\#})-L_1(Q_n)\}=O_P(n^{-1/2-\alpha/4})$. 
This proves that $d_{n1}(Q_n^{\#},Q_n)=O_P(n^{-1/2-\alpha/4})$. Using the same relation (\ref{handy}), this implies $d_{01}(Q_n^{\#},Q_n)=O_P(n^{-1/2-\alpha/4})$. 
$\Box$


\begin{lemma}\label{lemmadndo}
Suppose that $\int f^2_n dP_n=O_P(n^{-1/2-\alpha/4})$ and we know that $\pl f_n\pl_v^*<M$ for some $M<\infty$.
Then $\int f_n^2dP_0=O_P(n^{-1/2-\alpha/4})$.
\end{lemma}
{\bf Proof:}
We have
\begin{eqnarray*}
\int f_n^2 dP_0&=&-\int f_n^2d(P_n-P_0)+\int f_n^2 dP_n\\
&=&-\int f_n^2 d(P_n-P_0)+O_P(n^{-1/2-\alpha/4}).
\end{eqnarray*}
We have $\int f_n^2 d(P_n-P_0)=O_P(n^{-1/2})$.
This proves that $\int f_n^2 dP_0=O_P(n^{-1/2})$. By asymptotic equicontinuity of the empirical process indexed by cadlag functions with uniformly bounded sectional variation norm, it follows now also that $\int f_n^2 d(P_n-P_0)=O_P(n^{-1/2-\alpha/4})$ (the same proof can be found in Theorem 1 of \cite{vanderLaan15} using Lemma 10 in the same paper). Thus, this proves that indeed that $\int f_n^2 dP_0=O_P(n^{-1/2-\alpha/4})$ follows from $\int f_n^2dP_n=O_P(n^{-1/2-\alpha/4})$.
$\Box$

\section{Proof that the one-step TMLE $Q_n^{\#*}$ preserves rate of convergence of $Q_n^{\#}$}\label{AppendixC}

The following lemma establishes that the one-step TMLE $Q_n^{\#*}=Q_{n,\epsilon_n^{\#}}$ preserves the rate of convergence $d_{01}(Q_n^{\#},Q_0)=O_P(n^{-1/2-\alpha/4})$ of Theorem \ref{thnpbootmle} of $Q_n^{\#}$ in sense that also $d_{01}(Q_n^{\#*},Q_0)=O_P(n^{-1/2-\alpha/4})$. 
Recall the notation $L_1(Q_1,Q_2)=L_1(Q_1)-L_1(Q_2)$.
\begin{lemma}\label{lemmathalmleb}
Let $Q_n=\arg\min_{Q\in Q({\cal M})}P_n L_1(Q)$,  and let $\epsilon_n =\arg\min_{\epsilon}P_n L_1(Q_{n,\epsilon})$ for a parametric submodel $\{Q_{n,\epsilon}:\epsilon\}\subset {\cal M}$ thrugh $Q_n$ at $\epsilon =0$. 
Assume (\ref{sectionalvarbound}), (\ref{M20bounds}) so that we  know $d_{01}(Q_n,Q_0)=O_P(n^{-1/2-\alpha/4})$. 
By Lemma \ref{lemmathalmle} we also have $d_{01}(Q_n^*,Q_0)=O_P(n^{-1/2-\alpha/4})$. 
Let $Q_n^{\#}=\arg\min_{Q\in Q({\cal M})}P_n^{\#}L_1(Q)$ be the HAL-MLE on the bootstrap sample. 
By Theorem \ref{ThB} we also have $d_{n1}(Q_n^{\#},Q_n)=O_P(n^{-1/2-\alpha/4})$, where $d_{n1}(Q,Q_n)=P_n L_1(Q)-P_n L_1(Q_n)$, and 
$d_{01}(Q_n^{\#},Q_0)=O_P(n^{-1/2-\alpha/4})$.
Let $\epsilon_n^{\#}=\arg\min_{\epsilon} P_n^{\#}L_1(Q_{n,\epsilon}^{\#})$, and 
$Q_n^{\#*}=Q_{n,\epsilon_n^{\#}}^{\#*}$.

Then, 
\begin{equation}\label{thalmleb}
d_{01}(Q_n^{\#*},Q_0)=O_P(n^{-1/2-\alpha/4}).\end{equation}
\end{lemma}
{\bf Proof of Lemma \ref{lemmathalmleb}:}
Firstly, we note that
\begin{eqnarray*}
d_{01}(Q_n^{\#*},Q_0)&=&P_0L_1(Q_n^{\#*},Q_n^*)+P_0L_1(Q_n^*,Q_0)\\
&=&(P_0-P_n)L_1(Q_n^{\#*},Q_n^*)+P_n L_1(Q_n^{\#*},Q_n^*)+d_{01}(Q_n^*,Q_0)\\
&=&(P_0-P_n)L_1(Q_n^{\#*},Q_n^*)+d_{n1}(Q_n^{\#*},Q_n^*)+d_{01}(Q_n^*,Q_0)\\
&=&(P_0-P_n)L_1(Q_n^{\#*},Q_n^*)+d_{n1}(Q_n^{\#*},Q_n^*)+O_P(n^{-1/2-\alpha/4}).
\end{eqnarray*}
Using that $d_{n1}(Q_n^{\#},Q_n)$, $d_{01}(Q_n,Q_0)$, $d_{01}(Q_n^*,Q_0)$, and (thereby also, by \citep{vanderLaan15}) $(P_n-P_0)L_1(Q_n,Q_n^*)=O_P(n^{-1/2-\alpha/4})$ are all four $O_P(n^{-1/2-\alpha/4})$ we obtain
\begin{eqnarray*}
d_{n1}(Q_n^{\#*},Q_n^*)&=&P_nL_1(Q_{n,\epsilon_n^{\#}}^{\#} )-P_nL_1(Q_{n,\epsilon_n})\\
&=& P_nL_1(Q_{n,\epsilon_{n}^{\#}}^{\#},Q_n^{\#}) +P_nL_1(Q_n^{\#},Q_n)+P_nL_1(Q_n,Q_{n,\epsilon_n})\\
&=&(P_n-P_n^{\#})L_1(Q_{n,\epsilon_n^{\#}}^{\#},Q_n^{\#})
+P_n^{\#}L_1(Q_{n,\epsilon_n^{\#}}^{\#},Q_n^{\#})+d_{n1}(Q_n^{\#},Q_n)\\
&&+(P_n-P_0)L_1(Q_n,Q_n^*)+P_0L_1(Q_n,Q_0)+P_0L_1(Q_n^*,Q_0)\\
&\leq&(P_n-P_n^{\#})L_1(Q_{n,\epsilon_n^{\#}}^{\#},Q_n^{\#})+d_{n1}(Q_n^{\#},Q_n)+
(P_n-P_0)L_1(Q_n,Q_n^*)\\
&&+d_{01}(Q_n,Q_0)+d_{01}(Q_n^*,Q_0)\\
&=&(P_n-P_n^{\#})L_1(Q_{n,\epsilon_n^{\#}}^{\#},Q_n^{\#})+O_P(n^{-1/2-\alpha/4}).
\end{eqnarray*}
Plugging this bound for $d_{n1}(Q_n^{\#*},Q_n^*)$ in  our expression above for $d_{01}(Q_n^{\#*},Q_0)$ yields:
\begin{eqnarray*}
d_{01}(Q_n^{\#*},Q_0)&\leq& (P_0-P_n)L_1(Q_n^{\#*},Q_n^*)+(P_n-P_n^{\#})L_1(Q_n^{\#*},Q_n^{\#})+O_P(n^{-1/2-\alpha/4}).
\end{eqnarray*}
By assumption, we have that $L_1(Q_n^{\#*})$, $L_1(Q_n^{\#})$, $L_1(Q_n^*)$ are elements of the class of cadlag functions with universal bound on sectional variation norm, which is a uniform Donsker class. 
By empirical process theory for the empirical process $((P_n-P_0)f:f)$ and, conditional on $P_n$, for $((P_n^{\#}-P_n)f:f)$ indexed by this Donsker class, it follows that 
$d_{01}(Q_n^{\#*},Q_0)=O_P(n^{-1/2})+O_P(n^{-1/2-\alpha/4})$.
With this result in hand, we now revisit the 2 empirical process terms in the above bound for $d_{01}(Q_n^{\#*},Q_0)$ so that the $O_P(n^{-1/2})$ improves to $O_P(n^{-1/2-\alpha/4})$.
First, consider the second term. The $L^2(P_n)$-norm of $L_1(Q_n^{\#*},Q_n^{\#})$ is bounded by the sum of the $L^2(P_n)$-norms of $L_1(Q_n^{\#*},Q_0)$ and $L_1(Q_n^{\#},Q_0)$. The  $L^2(P_n)$-norm of $L_1(Q_n^{\#*},Q_0)$ is equivalent to $L^2(P_0)$-norm of $L_1(Q_n^{\#*},Q_0)$ (see Lemma \ref{lemmadndo}), which was just shown to be $O_P(n^{-1/4})$. 
The $L^2(P_n)$-norm of $L_1(Q_n^{\#},Q_0)$ can be bounded as sum of $L^2(P_n)$-norms of $L_1(Q_n^{\#},Q_n)$ and $L_1(Q_n,Q_0)$. These can be bounded in terms of $d_{n1}^{1/2}(Q_n^{\#},Q_n)$ and $d_{01}^{1/2}(Q_n,Q_0)$ (using Lemma \ref{lemmadndo} again). Thus, the $L^2(P_n)$ norm of $L_1(Q_n^{\#*},Q_n^{\#})$ is $O_P(n^{-1/4})$, so that we can establish again that $(P_n-P_n^{\#})L_1(Q_n^{\#*},Q_n^{\#})=O_P(n^{-1/2-\alpha/4})$.

Consider now the first empirical process term $(P_n-P_0)L_1(Q_n^{\#*},Q_n^*)$.
The $L^2(P_0)$-norm of $L_1(Q_n^{\#*},Q_n^*)$ can be bounded in terms of $L^2(P_0)$-norms of $L_1(Q_n^{\#*},Q_0)$ and $L_1(Q_n^*,Q_0)$, which thus is $O_P(n^{-1/4})$ as well. 
Therefore, it also follows that $(P_n-P_0)L_1(Q_n^{\#*},Q_n^*)=O_P(n^{-1/2-\alpha/4})$.
This proves that $d_{01}(Q_n^{\#*},Q_0)=O_P(n^{-1/2-\alpha/4})$.
$\Box$

\section{Proof of Theorem \ref{thnpboothaltmle}}
\label{AppendixD1}
Firstly, by definition of the remainder $R_{20}()$ we have the following two expansions:
\begin{eqnarray*}
\Psi(Q_n^{\#*})-\Psi(Q_0)&=&(P_n^{\#}-P_0) D^*(Q_n^{\#*},G_n^{\#})+R_{20}(Q_n^{\#*},G_n^{\#},Q_0,G_0)\\
&=& (P_n^{\#}-P_n)D^*(Q_n^{\#*},G_n^{\#})+(P_n-P_0)D^*(Q_n^{\#*},G_n^{\#})\\
&& +R_{20}(Q_n^{\#*},G_n^{\#},Q_0,G_0),\\
\Psi(Q_n^*)-\Psi(Q_0)&=&(P_n-P_0) D^*(Q_n^*,G_n)+R_{20}(Q_n^*,G_n,Q_0,G_0),
\end{eqnarray*}
where we ignored  $r_n=P_nD^*(Q_n^*,G_n)$ and its bootstrap analogue $r_n^{\#}=P_n^{\#}D^*(Q_n^{\#*},G_n^{\#})$ (which were both assumed to be $o_P(n^{-1/2})$).
Subtracting the first equality from the second equality yields:
\begin{eqnarray}
\Psi(Q_n^{\#*})-\Psi(Q_n^*) &=& (P_n^{\#}-P_n)D^*(Q_n^{\#*},G_n^{\#})\nonumber \\
&& + (P_n-P_0)\{D^*(Q_n^{\#*},G_n^{\#})-D^*(Q_n^*,G_n)\}\nonumber \\
&&+ R_{20}(Q_n^{\#*},G_n^{\#},Q_0,G_0)-R_{20}(Q_n^*,G_n,Q_0,G_0). \label{helpa}
\end{eqnarray}
Under the conditions of Theorem \ref{thefftmle}, we already established that $R_{20}(Q_n^*,G_n,Q_0,G_0)=O_P(n^{-1/2-\alpha/4})$.
By  assumption (\ref{boundingR2}), we can bound the first remainder $R_{20}(Q_n^{\#*},G_n^{\#},Q_0,G_0)$ by $f({\bf d}_{01}^{1/2}(Q_n^{\#*},Q_0),{\bf d}_{02}^{1/2}(G_n^{\#},G_0))$. 
Theorem \ref{thnpbootmle} established  that $d_{01}(Q_n^{\#},Q_0)=O_P(n^{-1/2-\alpha/4})$ and $d_{02}(G_n^{\#},G_0)=O_P(n^{-1/2-\alpha/4})$.  Using the fact that $f$ is a quadratic polyonomial, this now also establishes that
\begin{align*}
R_{20}(Q_n^{\#*},G_n^{\#},Q_0,G_0)=O_P(n^{-1/2-\alpha/4}).
\end{align*}


It remains to analyze the two leading empirical process terms in (\ref{helpa}). 
By our continuity assumption (\ref{contDstar}) on the efficient influence curve as function in $(Q,G)$, we have that  convergence of $d_{01}(Q_n^{\#*},Q_0)+d_{02}(G_n^{\#},G_0)$ to zero implies convergence of the square of the $L^2(P_0)$-norm of $D^*(Q_n^{\#*},G_n^{\#})-D^*(Q_0,G_0)$ at the same rate in probability. Since we already established convergence of $D^*(Q_n^*,G_n)-D^*(Q_0,G_0)$ to zero, this also establishes this result for the $L^2(P_0)$-norm of  $D^*(Q_n^{\#*},G_n^{\#})-D^*(Q_n^*,G_n)$.
By Lemma \ref{lemmadndo} this also proves that the $L^2(P_n)$-norm of the latter converges to zero in probability. 
By empirical process theory \citep{vanderVaart&Wellner11} (as in Appendix of \citep{vanderLaan15}), this teaches us that $(P_n^{\#}-P_n)D^*(Q_n^{\#*},G_n^{\#})=(P_n^{\#}-P_n)D^*(Q_n^*,G_n)+O_P(n^{-1/2-\alpha/4})$. This deals with the first leading term in (\ref{helpa}). 

By our continuity condition (\ref{contDstar}) we also have that 
\begin{align*}
P_0\{D^*(Q_n^{\#*},G_n^{\#})-D^*(Q_n^*,G_n)\}^2 \to_p 0    
\end{align*}
at this rate. Again, by \citep{vanderVaart&Wellner11} this shows 
$(P_n-P_0)\{D^*(Q_n^{\#*},G_n^{\#})-D^*(Q_n^*,G_n)\}=O_P(n^{-1/2-\alpha/4})$.
Thus we have shown  that
\[
\begin{array}{l}
(P_n^{\#}-P_n)D^*(Q_n^{\#*},G_n^{\#})+(P_n-P_0)\{D^*(Q_n^{\#*},G_n^{\#})-D^*(Q_n^*,G_n)\}\\
=(P_n^{\#}-P_n)D^*(Q_n^{*},G_n)
+O_P(n^{-1/2-\alpha/4}).\end{array}
\]

Thus, we have now shown, conditional on $(P_n:n\geq 1)$,
\[
n^{1/2}(\Psi(Q_n^{\#*})-\Psi(Q_n^*))=n^{1/2}(P_n^{\#}-P_n)D^*(Q_n^*,G_n)+o_P(1)\Rightarrow_d N(0,\sigma^2_0).
\]
This completes the proof of the Theorem for the HAL-TMLE. For a model ${\cal M}$ with extra structure (\ref{calFmodel}), this gives the result for the HAL-TMLE at the fixed $C^u$.  However, it follows straightforwardly that this proof applies uniformly in any $C$ in between $C_0$ and $C^u$, and thereby to a selector $C_n$ satisfying (\ref{Cn}).  $\Box$


\section{Understanding why $d_{n1}(Q_n^{\#},Q_n)$ is a quadratic dissimilarity}
\label{AppendixD}
\begin{lemma}\label{lemmadn1}
Assume extra model structure (\ref{calFmodel}) on ${\cal M}$.
Let $P_n R_{2L_1,n}(Q_n^{\#},Q_n)$ be defined as the exact second-order remainder of a first order Taylor expansion of $P_nL_1(Q)$ at $Q_n$:
\[ P_n \{L_1(Q_n^{\#})-L_1(Q_n)\}=P_n \frac{d}{dQ_n}L_1(Q_n)(Q_n^{\#}-Q_n)+P_n R_{2L_1,n}(Q_n^{\#},Q_n),\]
where $\frac{d}{dQ_n}L_1(Q_n)(h)=\left . \frac{d}{d\epsilon} L_1(Q_n+\epsilon h)\right |_{\epsilon =0}$ is the directional derivative in direction $h$.

We have $P_n\frac{d}{dQ_n}L_1(Q_n)(Q_n^{\#}-Q_n)\geq 0$ so that
\[
d_{n1}(Q_n^{\#},Q_n) \geq P_n R_{2L_1,n}(Q_n^{\#},Q_n).\]
\end{lemma}

In order to provide the reader a concrete example of what this empirical dissimilarity $d_{n1}(Q_n^{\#},Q_n)$ looks like, we provide here the corollary of Lemma \ref{lemmadn1} for the squared error loss.
\begin{corollary}
Consider the definitions of Lemma \ref{lemmadn1} and apply it to loss function
$L_1(Q)(O)=(Y-Q(X))^2$.
Then, $P_n R_{2L_1,n}(Q_n^{\#},Q_n)=P_n(Q_n^{\#}-Q_n)^2$, so that we have
\[
d_{n1}(Q_n^{\#},Q_n)\geq P_n (Q_n^{\#}-Q_n)^2 .\]
Since  $P_n  \{L_1(Q_n^{\#})-L_1(Q_n)\}^2=O_P( P_n (Q_n^{\#}-Q_n)^2)$, this implies
$P_n  \{L_1(Q_n^{\#})-L_1(Q_n)\}^2=O_P(d_{n1}(Q_n^{\#},Q_n))$.
\end{corollary}
{\bf Proof of Corollary:} 
We will prove $P_n R_{2L_1,n}(Q_n^{\#},Q_n)=P_n(Q_n^{\#}-Q_n)^2$. The remaining statement is then just an immediate corollary of Lemma \ref{lemmadn1}.
We have
\begin{eqnarray*}
d_{n1}(Q_n^{\#},Q_n)&=&\frac{1}{n}\sum_i \{2Y_iQ_n(X_i)-2Y_i Q_n^{\#}(X_i)+Q_n^{\#2}(X_i)-Q_n^2(X_i)\}\\
&=& \frac{1}{n}\sum_i\{2  (Q_n-Q_n^{\#})(X_i)  Y_i+Q_n^{\#2}(X_i)-Q_n^2(X_i)\}\\
&=&\frac{1}{n}\sum_i \{2(Q_n-Q_n^{\#})(X_i)(Y_i-Q_n(X_i))\\
&&\hfill +2(Q_n-Q_n^{\#})Q_n(X_i)+Q_n^{\#2}(X_i)-Q_n^2(X_i)\}\\
&=&\frac{1}{n}\sum_i 2(Q_n-Q_n^{\#})(X_i)(Y_i-Q_n(X_i))+\frac{1}{n}\sum_i (Q_n-Q_n^{\#})^2(X_i).\end{eqnarray*}
Note that the first term corresponds with $P_n\frac{d}{dQ_n}L_1(Q_n)(Q_n^{\#}-Q_n)$ and the second-order term with $P_nR_{2L_1,n}(Q_n^{\#},Q_n)$, where
$R_{2L_1,n}(Q_n^{\#},Q_n)=(Q_n^{\#}-Q_n)^2$. $\Box$

{\bf Proof of Lemma \ref{lemmadn1}:}
We need to prove that the linear approximation
\begin{align*}
P_n \frac{d}{dQ_n} L_1(Q_n)(Q_n^{\#}-Q_n) \leq 0.
\end{align*}
The extra model structure (\ref{calFmodel}) allows the explicit calculation of score equations for the HAL-MLE and its bootstrap analogue, which provides us then with the desired inequality.
 
Consider the  $h$-specific path  \[
Q_{n,\epsilon}^h(x)=(1+\epsilon h(0))Q_n(0)+\sum_s \int_{(0_s,x_s]}(1+\epsilon h_s(u_s)) dQ_{n,s}(u_s))\]
 for $\epsilon \in [0,\delta)$ for some $\delta>0$, where $h$ is uniformly bounded, and,  if $C^l<C^u$,  \[
 r(h,Q_n)\equiv  h(0)| Q_n(0)|+\sum_s \int_{(0_s,\tau_s]} h_s(u_s)| dQ_{n,s}(u_s)| \leq 0,\] while if $C^l=C^u$, then $r(h,Q_n)=0$.  Let ${\cal H}=\{h: r(h,Q_n)\leq 0,\pl h\pl_{\infty}<\infty\}$ be the set of possible functions $h$(i.e., functions of $s,u_s$), which defines a collection of paths $\{Q_{n,\epsilon}^h:\epsilon\}$ indexed by $h\in {\cal H}$. Consider a given $h\in {\cal H}$ and let's denote this path with $Q_{n,\epsilon}$, suppressing the dependence on $h$ in the notation.
For $\epsilon\geq 0$ small enough we have $(1+\epsilon h(0))>0$ and $1+\epsilon h_s(u_s)>0$. Thus, for $\epsilon\geq 0$ small enough we have
\begin{eqnarray*}
\pl Q_{n,\epsilon}\pl_v^*&=&(1+\epsilon h(0))| Q_n(0)|+\sum_s \int_{(0_s,\tau_s]}(1+\epsilon h_s(u_s))| dQ_{n,s}(u_s)| \\
&=& \pl Q_n\pl_v^*+\epsilon\left\{ h(0)| Q_n(0)|+\sum_s \int_{(0_s,\tau_s]}h_s(u_s) | dQ_{n,s}(u_s)| \right\} \\
&=&\pl Q_n\pl_v^*+\epsilon r(h,Q_n)\\
&\leq & \pl Q_n\pl_v^*,
\end{eqnarray*}
by assumption that $r(h,Q_n)\leq 0$. If $C^l=C^u$ and thus $r(h,Q_n)=0$, then the above shows $\pl Q_{n,\epsilon}\pl_v^*=\pl Q_n\pl_v^*$.
Thus, for a small enough $\delta>0$ $\{Q_{n,\epsilon}:0\leq \epsilon<\delta\}$ represents a path of cadlag functions with sectional variation norm bounded from below and above: $C^l \leq \pl Q_n\pl_v^*\leq C^u$. In addition, we have that $dQ_{n,s}(u_s)=0$ implies $(1+\epsilon h_s(u_s))dQ_{n,s}(u_s)=0$ so that the support of $Q_{n,\epsilon}$ is included in the support $A$ of $Q_n$ as defined by ${\cal F}_A^{np}$. Thus, this proves that for $\delta>0$ small enough this path $\{Q_{n,\epsilon}:0\leq \epsilon\leq \delta\}$ is indeed a  submodel of the parameter space of $Q$, defined as ${\cal F}_{A}^{np}$ or ${\cal F}_A^{np+}$.

We also have that
 \[
 Q_{n,\epsilon}-Q_n=\epsilon\left\{ Q_n(0)h(0)+ \sum_s \int_{(0_s,x_s]} h_s(u_s)dQ_{n,s}(u_s)\right\} .\]
 Thus, this path generates a direction $f(h,Q_n)$  at $\epsilon=0$ given by:
 \[
 \frac{d}{d\epsilon}Q_{n,\epsilon}=f(h,Q_n)\equiv Q_n(0)h(0)+ \sum_s \int_{(0_s,x_s]} h_s(u_s)dQ_{n,s}(u_s) .\]
 Let ${\cal S}\equiv \{f(h,Q_n): h\in {\cal H}\}$ be the collection of directions generated by our family of paths.
 By definition of the MLE $Q_n$, we also have that $\epsilon \rightarrow P_n L_1(Q_{n,\epsilon})$ is minimal over $[0,\delta)$ at $\epsilon =0$.
This shows that the derivative of $P_n L_1(Q_{n,\epsilon})$ from the right at $\epsilon =0$ is  non-negative:
 \[
 \frac{d}{d\epsilon^+ }P_n L_1(Q_{n,\epsilon })\geq 0\mbox{ at $\epsilon =0$}.\]
 This derivative is given by $P_n \frac{d}{dQ_n}L_1(Q_n)(f(h,Q_n))$, where $d/dQ_nL_1(Q_n)(f(h,Q_n))$ is the directional (Gateaux) derivative of $Q\rightarrow L_1(Q)$ at $Q_n$ in  in direction $f(h,Q_n)$.
  Thus for each $h\in {\cal H}$, we have
\[
P_n \frac{d}{dQ_n}L_1(Q_n)(f(h,Q_n)) \geq 0 .\]
Suppose that
\begin{equation}\label{keya1}
Q_n^{\#}-Q_n\in {\cal S}=\{f(h,Q_n):h\in {\cal H}\}.\end{equation}
Then, we have
\[
P_n \frac{d}{dQ_n}L_1(Q_n)(Q_n^{\#}-Q_n)\geq 0.\]
Combined with the stated second-order Taylor expansion of $P_n L_1(Q)$ at $Q=Q_n$ with exact second-order remainder  $P_nR_{2L_1,n}(Q_n^{\#},Q_n)$, this proves
\[
P_n\{L_1(Q_n^{\#})-L_1(Q_n)\}\geq P_nR_{2L_1,n}(Q_n^{\#},Q_n).\]
Thus it remains to show (\ref{keya1}).

In order to prove (\ref{keya1}), let's solve explicitly for $h$ so that $Q_n^{\#}-Q_n=f(h,Q_n)$ and then verify that $h\in {\cal H}$ satisfies its assumed constraints (i.e., $r(h,Q_n)\leq 0$ if $C^l<C^u$ or $r(h,Q_n)=0$ if $C^l=C^u$,  and $h$ is uniformly bounded).
We have
\begin{eqnarray*}
Q_n^{\#}-Q_n&=&Q_n^{\#}(0)-Q_n(0)+\sum_s\int_{(0_s,x_s]} d(Q_{n,s}^{\#}-dQ_{n,s})(u_s)\\
&=& Q_n^{\#}(0)-Q_n(0)+\sum_s \int_{(0_s,x_s]} \frac{d(Q_{n,s}^{\#}-dQ_{n,s})}{dQ_{n,s}} dQ_{n,s}(u_s),
\end{eqnarray*}
where we used that $Q_{n,s}^{\#}\ll Q_{n,s}$ for each subset $s$.
Let $h(Q_n^{\#},Q_n)$ be defined by
\begin{eqnarray*}
h(Q_n^{\#},Q_n)(0)&=&(Q_n^{\#}(0)-Q_n(0))/Q_n(0)\\
h_s(Q_n^{\#},Q_n)&=&\frac{d(Q_{n,s}^{\#}-dQ_{n,s})}{dQ_{n,s}}\mbox{ for all subsets $s$} .
\end{eqnarray*}
For this choice $h(Q_n^{\#},Q_n)$, we have $f(h,Q_n)=Q_n^{\#}-Q_n$.
First, consider the case $Q({\cal M})={\cal F}_A^{np}$ or $Q({\cal M})={\cal F}_A^{np+}$, but $C^l<C^u$.
We now need to verify that $r(h,Q_n)\leq 0$ for this choice $h=h(Q_n^{\#},Q_n)$.
We have
\begin{eqnarray*}
r(h,Q_n)&=&\frac{Q_n^{\#}(0)-Q_n(0)}{Q_n(0)}| Q_n(0)|+\sum_s\int_{(0_s,\tau_s]}\frac{dQ_{n,s}^{\#}-dQ_{n,s}}{dQ_{n,s}}| dQ_{n,s}|\\
&=&I(Q_n(0)>0)\{Q_n^{\#}(0)-Q_n(0)\}+I(Q_n(0)\leq 0)\{Q_n(0)-Q_n^{\#}(0)\}\\
&&+\sum_s\int_{(0_s,\tau_s]}I(dQ_{n,s}\geq 0) d(Q_{n,s}^{\#}-dQ_{n,s})\\
&&+\sum_s\int_{(0_s,\tau_s]} I(dQ_{n,s}<0)d(Q_{n,s}-Q_{n,s}^{\#})\\
&=&-\pl Q_n\pl_v^*+Q_n^{\#}(0)\{ I(Q_n(0)>0)-I(Q_n(0)\leq 0)\} \\
&&+\sum_s \int_{(0_s,\tau_s]}\{I(dQ_{n,s}\geq 0) -I(dQ_{n,s}\leq 0\} dQ_{n,s}^{\#}\\
&\leq&-\pl Q_n\pl_v^*+| Q_n^{\#}(0)| +\sum_s\int_{(0_s,\tau_s]} | dQ_{n,s}^{\#}(u_s)| \\
&=&-\pl Q_n\pl_v^*+\pl Q_n^{\#}\pl_v^*\\
&\leq &0,
\end{eqnarray*}
since $\pl Q_n^{\#}\pl_v^*\leq \pl Q_n\pl_v^*$, by assumption.
Thus, this proves that indeed $r(h,Q_n)\leq 0$ and thus that $Q_n^{\#}-Q_n\in {\cal S}$.
Consider now the case that $Q({\cal M})={\cal F}_A^{np+}$ {\em and } $C^l=C^u$. Then $\pl Q_n\pl_v^*=\pl Q_n^{\#}\pl_v^*=C^u$. We now need to show that $r(h,Q_n)=0$ for this choice $h=h(Q_n^{\#},Q_n)$. We now use the same three equalities as above, but now use that $dQ_{n,s}(u_s)\geq 0$ and $Q_n(0)\geq 0$, by definition of ${\cal F}_A^{np+}$, which then shows $r(h,Q_n)=0$.

This proves (\ref{keya1}) and thereby  completes the proof of Lemma \ref{lemmadn1}. $\Box$

\section{Number of Non-zero HAL Coefficients As a Function of Sample Size}
\label{AppendixF}
\begin{figure}[ht]
  \centering
  \includegraphics[width=0.45\textwidth]{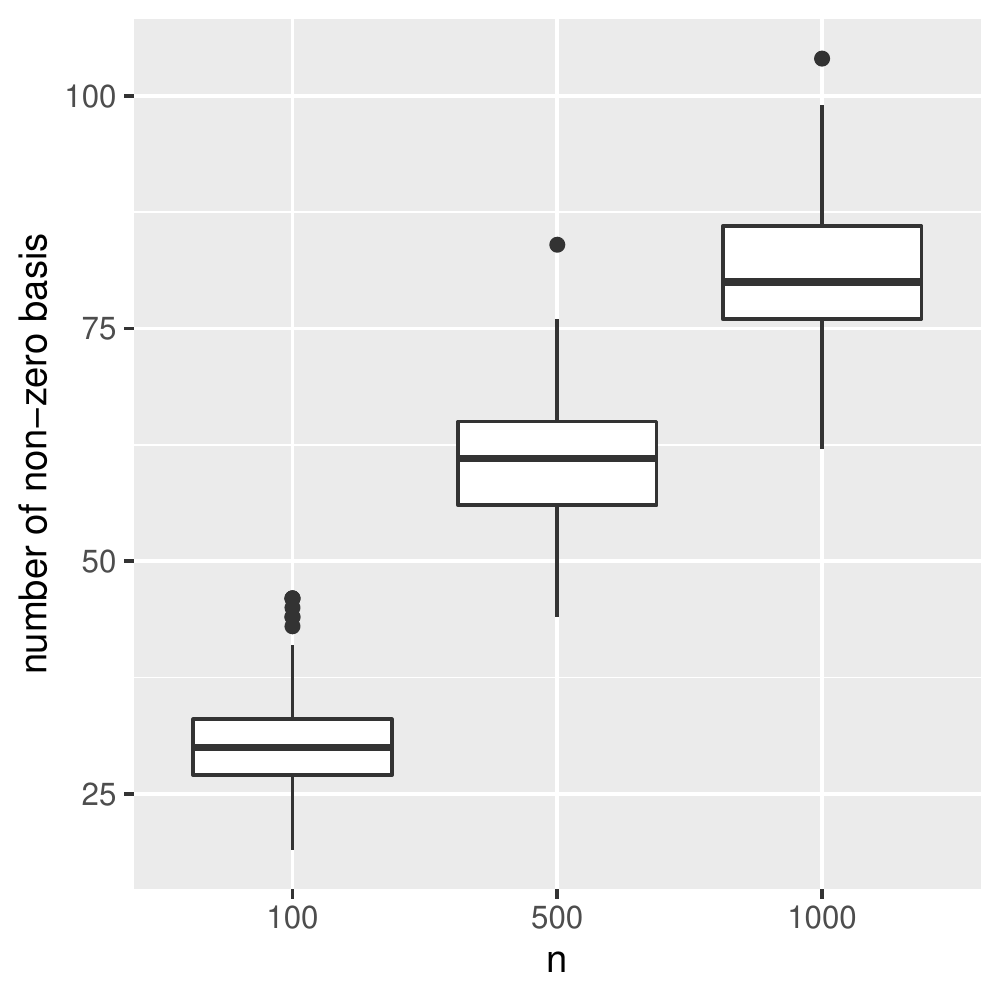}
  \includegraphics[width=0.45\textwidth]{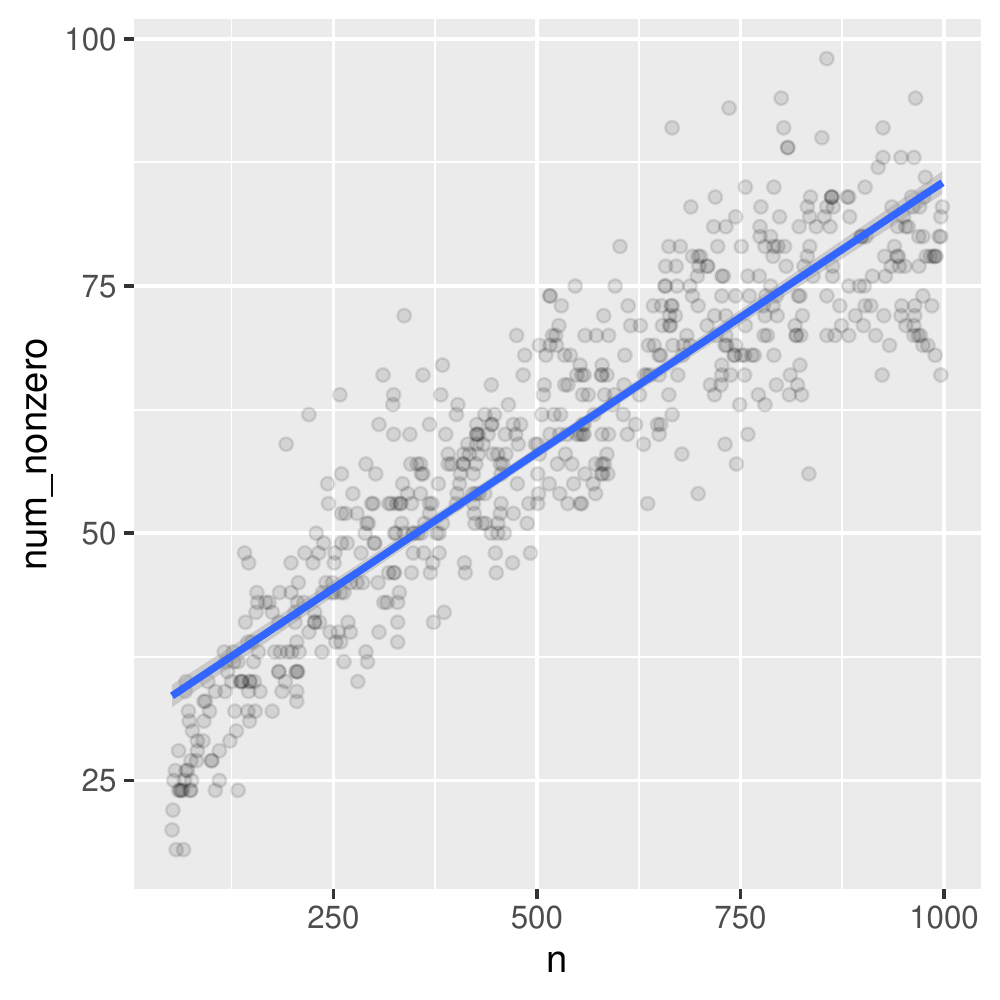}
  \caption{The number of non-zero coefficients in $Q_n$ as a function of sample size (under simulation 1 at $a_1 = 0.5$)}
\end{figure}

\end{document}